\documentclass[reqno,letterpaper,11pt]{amsart}
\usepackage[letterpaper, margin=1in]{geometry}

\usepackage[utf8]{inputenc}
\usepackage[T1]{fontenc}
\usepackage[english]{babel}
\usepackage{csquotes}
\usepackage[babel]{microtype}
\usepackage[dvipsnames]{xcolor}

\usepackage{amsmath}
\usepackage{amssymb}
\usepackage{amsthm}
\usepackage{mathrsfs}
\usepackage{mathtools}
\usepackage{commath}
\usepackage{thmtools}

\usepackage{bbm}
\usepackage{dsfont}
\usepackage{lmodern}
\usepackage{fixcmex}
\usepackage{cases}
\usepackage{centernot}
\usepackage{yhmath}

\usepackage{enumitem}
\setlist[enumerate,1]{label*=\arabic*., leftmargin=2em}

\usepackage[pdfborderstyle={/S/U/W 0},unicode]{hyperref}
\hypersetup{
    colorlinks=true,
    linkcolor=blue!85!black,
    citecolor=blue!85!black,
    urlcolor=blue!85!black,
}

\usepackage{cleveref}

\numberwithin{equation}{section}
\makeatletter
\@ifundefined{newcounteralias}{}{%
  \renewcommand\thmt@autorefsetup{\@xa\def\csname\thmt@envname autorefname\@xa\endcsname\@xa{\thmt@thmname}}%
}
\makeatother

\declaretheoremstyle[
  headfont=\normalfont\bfseries,
  bodyfont=\normalfont,
]{noital}

\declaretheorem[style=plain,name=Theorem,numberwithin=section]{theorem}
\declaretheorem[style=plain,sibling=theorem,name=Proposition]{proposition}
\declaretheorem[style=plain,sibling=theorem,name=Lemma]{lemma}
\declaretheorem[style=plain,sibling=theorem,name=Corollary]{corollary}

\declaretheorem[style=noital,sibling=theorem,name=Remark]{remark}
\declaretheorem[style=noital,sibling=theorem,name=Definition]{definition}

\newcommand{\defined}{\mathrel{\coloneqq}}
\DeclarePairedDelimiter{\p}{\lparen}{\rparen}

\renewcommand{\leq}{\leqslant}

\renewcommand{\geq}{\geqslant}
\newcommand{\st}{\mathbin{\colon}}
\DeclarePairedDelimiter{\card}{\lvert}{\rvert}
\DeclareMathOperator{\ind}{\mathbf{1}}
\newcommand{\from}{\colon}
\DeclareMathOperator{\compo}{\circ}
\DeclareMathOperator{\Aut}{Aut}
\DeclarePairedDelimiter{\floor}{\lfloor}{\rfloor}
\DeclarePairedDelimiter{\ceil}{\lceil}{\rceil}
\newcommand{\isomorphic}{\simeq}
\newcommand{\fourier}{\widehat}
\DeclareMathOperator{\loglog}{log\,log}
\DeclareMathOperator{\logloglog}{log\,log\,log}
\DeclareMathOperator{\lcm}{lcm}
\newcommand{\divides}{\mathrel{\mid}}

\DeclareMathOperator{\eqmod}{\equiv}
\newcommand*{\conj}[1]{\overline{#1}}
\DeclareMathOperator{\conv}{\ast}
\DeclareMathOperator{\supp}{supp}

\undef{\set}
\DeclarePairedDelimiter{\set}{\lbrace}{\rbrace}

\undef{\emptyset}
\newcommand{\emptyset}{\varnothing}

\undef{\mod}
\newcommand{\mod}[1]{\ (\mathrm{mod}\ #1)}

\undef{\abs}
\DeclarePairedDelimiterX{\abs}[1]
  {\lvert}{\rvert}{\ifblank{#1}{\,\cdot\,}{#1}}

\undef{\norm}
\DeclarePairedDelimiterX{\norm}[1]
  {\lVert}{\rVert}{\ifblank{#1}{\,\cdot\,}{#1}}

\DeclarePairedDelimiterX{\inner}[2]
  {\langle}{\rangle}{\ifblank{#1}{\,\cdot\,}{#1},\ifblank{#2}{\,\cdot\,}{#2}}

\newcommand{\eps}{\varepsilon}

\newcommand{\CC}{\mathbb{C}}
\newcommand{\FF}{\mathbb{F}}
\newcommand{\NN}{\mathbb{N}}
\newcommand{\RR}{\mathbb{R}}
\let\SS\relax
\newcommand{\SS}{\mathbb{S}}
\newcommand{\TT}{\mathbb{T}}
\newcommand{\ZZ}{\mathbb{Z}}

\newcommand{\cB}{\mathcal{B}}
\newcommand{\cC}{\mathcal{C}}
\newcommand{\cF}{\mathcal{F}}
\newcommand{\cI}{\mathcal{I}}
\newcommand{\cS}{\mathcal{S}}

\newcommand{\rL}{\mathrm{L}}
\newcommand{\rI}{\mathrm{I}}
\newcommand{\rK}{\mathrm{K}}

\usepackage[backend=biber, style=numeric, maxnames=99, giveninits=true, doi=false, isbn=false, url=false]{biblatex}
\bibliography{references-arxiv}

\newcommand{\wiener}{{\texttt{A}}}
\DeclareMathOperator{\str}{str}
\newcommand{\dilation}{\cdot}
\newcommand{\Lin}{\mathsf{Lin}}
\DeclareMathOperator{\Gr}{Gr}
\DeclareMathOperator{\Tr}{Tr}

\begin{document}

\title{Asymptotic independence of solution-free constraints over finite fields}

\author{Arker Oke Soe}
\author{Victor Souza}
\address{Department of Mathematics, Cornell University, Ithaca NY, USA}
\email{as3738@cornell.edu,vsouza@cornell.edu}
\address{Sidney Sussex College, University of Cambridge, Cambridge, UK}
\email{vss28@cam.ac.uk}

\begin{abstract}
If $L$ and $L'$ are two linear forms in least three variables and integer coefficients, we show that the constraints `\emph{$A$ has no solution to $L(x) = 0$}' and `\emph{$A^{-1}$ has no solution to $L'(x) = 0$}', are asymptotically independent for subsets $A$ of a finite field.
More precisely, if $\alpha$ and $\beta$ are respectively the largest densities of subsets of $\mathbb{F}_q$ with no solutions to $L(x) = 0$ and $L'(x) = 0$, then every set satisfying both the conditions mentioned above has density at most $\alpha \beta + o(1)$ and the number of such set is $2^{(\alpha \beta + o(1))q}$.
We show a similar independence result for diagonal equations of the form $a_1 x_1^r + \dotsb + a_s x_s^r = 0$ and even multiplicative constraints such as `\emph{$A$ is product-free}'.
\end{abstract}

\maketitle


\section{Introduction}
 
The number of primes in the interval $[1,x]$ is asymptotically $x/\log x$, while the the number such primes that are congruent to $a$ modulo $b$ is asymptotically $x/(\varphi(b) \log x)$ if $a$ and $b$ are coprime.
In other words, the condition that `\emph{$n$ is prime}' and the condition that `\emph{$n \eqmod a \mod{b}$}' are asymptotically independent in the set of integers coprime to $b$.
Several other results in number theory can be cast in a similar way.
To name a few, the Chinese remainder theorem says that `$n \eqmod a \mod{b}$' and `$n \eqmod c \mod{d}$' are independent if $b$ and $d$ are coprime.
The equidistribution of the sequence $\set{\alpha n}_{n \geq 0}$ in $\RR/\ZZ$ when $\alpha$ is irrational says that the conditions `\emph{$\set{\alpha n} \in I$}' and `\emph{$x \in I$}' are asymptotically independent for every measurable $I \subseteq \RR/\ZZ$.

In this work, we investigate the asymptotic independence of constraints of the form `\emph{$A$ has no solution to a certain equation}' when $A$ is a subset of a finite field.


\subsection{Sets without solutions to linear equations}

Let $L(x_1,\dotsc,x_s) = a_1 x_1 + \dotsb + a_s x_s$ be a linear form with non-zero integer coefficients.
We say that a subset $A$ of a commutative group $G$ is $L$-free if there is no $(x_1, \dotsc, x_s) \in A^s$ with $L(x_1, \dotsc, x_s) = 0$.
For $S \subseteq G$, we write
\begin{equation*}
    \cF_L(S) \defined \set[\big]{A \subseteq S \st \text{$A \subseteq S$ is $L$-free }}
\end{equation*}
for the family of $L$-free subsets of $S$, and for finite subsets $S \subseteq G$, we write
\begin{equation*}
    d_L(S) \defined \max\set[\big]{\, \card{A}/\card{S} \st \text{$A \subseteq S$ is $L$-free }}
\end{equation*}
for the maximum density of an $L$-free subset of $S$.
We also denote by
\begin{equation*}
    m_L(S) \defined \card[\big]{\set{ (x_1, \dotsc, x_s) \in S^s \st L(x_1, \dotsc, x_s) = 0}} / \card{G}^s
\end{equation*}
for the normalised count of solutions to $L(x) = 0$ in a set $S \subseteq G$.
If $G$ is a compact commutative group equipped with the Haar probability measure $\mu_G$, we define $d_L(G) \defined \sup\set{ \mu_G(A) \st \text{$A$ is $L$-free }}$.

While the systematic study of $d_L(S)$ was initiated by Ruzsa~\cite{Ruzsa1993-py,Ruzsa1995-lv}, many classical results in additive combinatorics, such as Roth's theorem~\cite{Roth1953-nk}, the study of Sidon sets~\cite{Erdos1941-ru} and sum-free sets~\cite{Schur1916-lc,Erdos1965-on}, can be cast in this framework.

We write $\Lin_s$ for the set of linear forms $L(x_1,\dotsc,x_s) = a_1 x_1 + \dotsb + a_s x_s$ with non-zero integer coefficients and $s$ variables.
Write $\Lin_{\geq s} \defined \bigcup_{k \geq s} \Lin_s$ and $\Lin \defined \Lin_{\geq 2}$.
In our work, we will be mainly concerned with finite fields $\FF_q$, where $q = p^d$ for some prime $p$ and $d \geq 1$.
We say that a for $L \in \Lin$ is \emph{degenerate in $\FF_q$} if some of its coefficients are equal to zero in $\FF_q$, and \emph{non-degenerate\footnote{The linear form $x_1 - 3 x_2 + 7 x_3$ is degenerate in $\FF_{7^d}$ as it should be treated as a form in two variables.} in $\FF_q$} otherwise.
The linear form $L \in \Lin$ is \emph{translation-invariant in $\FF_q$} whenever $a_1 + \dotsb + a_s = 0$ in $\FF_q$, and in this case, only the empty set is $L$-free since $L(x, \dotsc, x) = 0$ for every $x \in \FF_q$.
It is customary to classify such solutions as \emph{trivial} and to deem a set $L$-free if it has no non-trivial solutions.
Even then, it follows from density increment arguments that $L$-free sets are sparse: that is, $d_L(\FF_q) = o(1)$ as $q \to \infty$.
We only  consider linear forms that are not translation-invariant, and many of our results become vacuous otherwise.

The additive group of $\FF_q$ is isomorphic to $\FF_p^d$, so $d_L(\FF_{p^d}) = d_L(\FF_p^d)$.
Whenever $L \in \Lin_{\geq 3}$ is not translation-invariant, the value of $d_L(\FF_q)$ does not converge to $0$ when $q \to \infty$, however, to make sense of the possible asymptotic values, it is sensible to consider two different scenarios.
In the \emph{fixed dimension regime}, we consider $q = p^d$ when $p \to \infty$ and $d$ is fixed, while in the \emph{fixed characteristic regime}, we consider $q = p^d$ when $p$ is fixed and $d \to \infty$.
In any case, $d_L(\FF_q)$ converges to some positive value.
Candela and Sisask~\cite{Candela2011-cg} have shown that
\begin{equation}
\label{eq:density-L-Fp}
    \lim_{p \to \infty} d_L(\FF_p) = d_L(\TT)
\end{equation}
for every $L \in \Lin_{\geq 3}$, where $\TT \defined \RR/\ZZ$ is the torus equipped with the Lebesgue probability measure.
It is not hard (see \Cref{cor:candela-sisask-d}) to then conclude that $\lim_{p \to \infty} d_L(\FF_p^d) = d_L(\TT^d) = d_L(\TT)$ for every fixed $d \geq 1$.
It is also established in \Cref{prop:limit-fix-char} that, for every $L \in \Lin_{\geq 3}(\FF_p)$, 
\begin{equation}
\label{eq:density-L-Fq}
    \lim_{d \to \infty} d_L(\FF_p^d) = d_L(\FF_p^\infty),
\end{equation}
where $\FF_p^\infty$ is the direct product of countably many copies of $\FF_p$ equipped with the product probability Haar measure and product topology.
While we can show that the limits \eqref{eq:density-L-Fp} and \eqref{eq:density-L-Fq} exists, determining their value is a completely different matter.

The study of $L$-free sets for the linear form $\Sigma(x,y,z) \defined x + y - z$ goes back to the work of Schur~\cite{Schur1916-lc}.
A subset that is $\Sigma$-free is called \emph{sum-free}, so we denote by $d_{\Sigma}(S)$ for the density of the largest sum-free subset of $S$.
Green and Ruzsa~\cite{Green2005-no} completely determined\footnote{Building of previous work of Diananda and Yap~\cite{Diananda1969-ia}, Rhemtulla and Street~\cite{Rhemtulla1970-px}, and of Yap~\cite{Yap1972-rb, Yap1975-vb}.} the densities $d_\Sigma(G)$ for all finite commutative groups $G$.
For finite fields, we have $d_\Sigma(\FF_p^d) = \floor{(p+1)/3}/p$ and thus $d_\Sigma(\FF_p^d) = 1/3 + o(1)$ as $p \to \infty$.
Moreover, the maximum density is attained by sets of the form $I_p \times \FF_p^{d-1}$ where $I_p$ is an appropriate \emph{middle-third} interval with $\card{I_p} = \floor{(p+1)/3}$.
In the same work, Green and Ruzsa also determined that $\log_2 \card{\cF_\Sigma(G)} = (d_\Sigma(G) + o(1))\card{G}$, which can been seen as a weak form of the Cameron-Erdős conjecture\footnote{Cameron and Erdős~\cite{Cameron1990-ez} asked, among many other questions, for an asymptotic enumeration of $\cF_\Sigma([N])$ where $[N] \defined \set{1, \dotsc, N} \subseteq \ZZ$. This was accomplished independently by Green~\cite{Green2004-tx} and Sapozhenko~\cite{Sapozhenko2003-zs,Sapozhenko2008-fq}, who showed $\card{\cF_\Sigma([N])} = (C_{N} + o(1))2^{N/2}$ where $C_N$ is a constant that depends solely on the parity of $N$. Calkin~\cite{Calkin1990-zb}, Alon~\cite{Alon1991-ab}, and Erdős and Granville had previously obtained the weaker estimate $\card{\cF_\Sigma([N])} = 2^{(1/2 + o(1))N}$. Improving on results of Lev and Schoen~\cite{Lev2002-pp} in the prime field case, and Green and Ruzsa~\cite{Green2004-yr} showed that $\card{\cF_\Sigma(\FF_p)} = 2^{(1/3 + o(1))p}$. Green and Ruzsa~\cite{Green2005-no} asymptotically determined $\card{\cF_\Sigma(G)}$ if $G$ is of \emph{Type $\rI(p)$}, generalising a result of Lev, Łuczak, and Schoen~\cite{Lev2001-pj} and of Sapozhenko~\cite{Sapozhenko2002-xc}, who did the same when $G$ has even order.}.

One can also establish the counting result more generally.
Indeed, for $L \in \Lin(\FF_q)_{\geq 3}$, we have
\begin{equation}
\label{eq:card-family-L-Fq}
    \log_2 \card{\cF_L(\FF_q)} = (d_{L}(\FF_q) + o(1))q
\end{equation}
as $q \to \infty$.
This can be found, for instance, in the work of  Green~\cite{Green2005-nn} via his arithmetical removal lemma, or in Saxton and Thomason~\cite{Saxton2015-cr,Saxton2016-xf} via containers.
Note that in \eqref{eq:density-L-Fp}, \eqref{eq:density-L-Fq}, and in \eqref{eq:card-family-L-Fq} we assume that $L \in \Lin(\FF_q)_{\geq 3}$.
In \Cref{prop:two-variables} we show that these limits often do not exist for $L \in \Lin(\FF_q)_2$, an even when they do, \eqref{eq:card-family-L-Fq} does not hold.


\subsection{Self-inverse sets with linear constraints}

We study now the effect of multiplicative structure on large $L$-free sets.
The first problem we consider is the following: what is the maximum density of an self-inverse $L$-free subset $A \subseteq \FF_q$?
More generally, let $R$ be a finite ring and denote $R^\times$ for the set of units (invertible elements) in $R$.
For a set $A \subseteq R^\times$, denote by $A^{-1} \defined \set{x^{-1} \st x \in A}$, and we say that $A \subseteq R$ is \emph{self-inverse} if $A = A^{-1}$.
For a finite subset $S \subseteq R$, we are then concerned with the following quantity
\begin{equation*}
    d_L^{\times}(S) \defined \max\set[\big]{ \card{A}/\card{S} \st \text{$A \subseteq S$ is $L$-free and $A = A^{-1}$}}.
\end{equation*}
This quantity was considered first in the sum-free case by Bienvenu, Hennecart, and Shkredov~\cite{Bienvenu2019-sz}, who showed that $1/9 - o(1) \leq d_\Sigma^\times(\FF_p) \leq 0.3051 + o(1)$ and as posed the question of determining $d_\Sigma^\times(\FF_p)$.
Benjamin~\cite{Benjamin2021-cq} then showed that $d_\Sigma^\times(\FF_p) \leq 1/4 - \eps + o(1)$ for some fixed $\eps >0$ and showed that $d_\Sigma^\times(\FF_{2^d}) = 1/4 + o(1)$, concluding then that $d_\Sigma^\times(\FF_q)$ does not converge as $q \to \infty$.
However, she conjectured that $d_\Sigma^\times(\FF_q) = d_\Sigma(\FF_q)^2 + o(1)$ as $q \to \infty$.
Semchankau~\cite{Semchankau2025-nt} determined that $d_\Sigma^\times(\FF_{p}) = 1/9 + o(1)$, which settles the conjecture in the prime field case, while also introducing the technique of \emph{wrappers}, which may be though as a refinement of the \emph{granularisation} method of Green and Ruzsa~\cite{Green2004-yr}.
In our first result, we determine $d_L^\times(\FF_q)$ asymptotically for every linear form $L$ in at least three variables.

\begin{theorem}
\label{thm:density-equation-inverse-fq}
For every $L \in \Lin_{\geq 3}$, we have that
\begin{equation*}
    d_L^\times(\FF_q) = d_L(\FF_q)^2 + o(1)
\end{equation*}
when $q \to \infty$, as long as $L$ is non-degenerate in $\FF_q$.
\end{theorem}

In particular, this gives
\begin{equation*}
    d_L^\times(\FF_{p^d}) = \begin{cases}
        d_L(\TT)^2 + o(1) & \text{if $d$ is fixed and $p \to \infty$,} \\
        d_L(\FF_p^\infty)^2 + o(1) & \text{if $p$ is fixed and $d \to \infty$.}
    \end{cases}
\end{equation*}
Benjamin's conjecture is then settled in full, and we have $d_\Sigma^\times(\FF_{p^d}) = 1/9 + o(1)$ as $p \to \infty$ for every fixed $d \geq 1$, and $d_\Sigma^\times(\FF_{p^d}) = \floor{(p + 1)/3}^2/p^2 + o(1)$ as $d \to \infty$ for every fixed $p \geq 2$.

\Cref{thm:density-equation-inverse-fq} can be interpreted as follows: the condition that `\emph{$A$ is sum-free}', and the condition that `\emph{$A^{-1}$ is sum-free}' are asymptotically \emph{independent} and as each incur in a density cost of $d_L(\FF_q)$, the total density cost is $d_L(\FF_q)^2$.
This suggest that the condition that $A = A^{-1}$ is not fully necessary.
Indeed, we establish a more general independence phenomena.
For $S \subseteq R$, define
\begin{equation*}
    d_{L,L'}(S) \defined \max\set[\big]{ \card{A}/\card{S} \st \text{$A \subseteq S \cap R^\times$ is $L$-free and $A^{-1}$ is $L'$-free }},
\end{equation*}
and also define the family
\begin{equation*}
    \cF_{L,L'}(S) \defined \set[\big]{ A \subseteq S \cap R^\times  \st \text{$A$ is $L$-free and $A^{-1}$ is $L'$-free }}.
\end{equation*}
We have the following independence result.

\begin{theorem}
\label{thm:density-independence-fq}
For every $L, L' \in \Lin(\FF_q)_{\geq 3}$, we have that
\begin{equation*}
    d_{L,L'}(\FF_q) = d_{L}(\FF_q) \, d_{L'}(\FF_q) + o(1)
\end{equation*}
when $q \to \infty$, as long as $L$ and $L'$ are non-degenerate in $\FF_q$.
\end{theorem}

While \Cref{thm:density-independence-fq} and \Cref{thm:density-equation-inverse-fq} indicate that self-inverse sum-free sets are related to sum-free sets whose inverse set is sum-free, the family $\cF_{L,L}(\FF_q)$ is much larger than $\cF_L^\times(\FF_q)$.
Indeed, we have the following enumeration results.

\begin{theorem}
\label{thm:counting-independence-fq}
For every $L, L' \in \Lin(\FF_q)_{\geq 3}$, we have that
\begin{equation}
\label{eq:counting-independence-fq}
    \log_2 \card{\cF_{L,L'}(\FF_q)} = \p[\big]{d_{L}(\FF_q) \, d_{L'}(\FF_q) + o(1)}q,
\end{equation}
when $q \to \infty$, as long as $L$ and $L'$ are non-degenerate in $\FF_q$.
\end{theorem}

The lower bound in \eqref{eq:counting-independence-fq} is simple: just consider all $2^{\card{X}}$ subsets of a set $X \in \cF_{L,L'}(\FF_q)$ with maximum cardinality.
The family $\cF_L^\times(\FF_q)$ is indeed smaller than $\cF_{L,L}(\FF_q)$.

\begin{theorem}
\label{thm:counting-equation-inverse-fq}
For every $L \in \Lin(\FF_q)_{\geq 3}$, we have that
\begin{equation*}
    \log_2 \card{\cF_L^\times(\FF_q)} = \p[\big]{d_L^\times(\FF_q)/2 + o(1)} q,
\end{equation*}
when $q \to \infty$, as long as $L$ is non-degenerate in $\FF_q$.
\end{theorem}

The factor of two comes from the fact that element of a self-inverse subset of $\FF_q^\times$ come in pairs $\set{x,x^{-1}}$, leading to half as many choices of subsets.


\subsection{Asymptotic independence for diagonal constraints}

For our grand independence result, we consider the following definitions.
For an integer $r \neq 0$ and a subset $A \subseteq R$ of a ring, define\footnote{This notation is not standard and we emphasise that $A^2$ does not denote the product set $A \cdot A \defined \set{x y \st x , y \in A}$.} the set of $r$-powers of elements in $A$ by $A^r \defined \set{a^r \st a \in A}$.
If $r < 0$, then we assume $A$ is a subset of the units $R^\times$.
Also denote the dilation as $r \dilation A \defined \set{r a \st a \in A}$.
For $L_1, \dotsc, L_t \in \Lin$ and $r_1, \dotsc, r_t \in \ZZ$, and $S \subseteq G$, define the family
\begin{equation*}
    \cF_{L_1, \dotsc, L_t}^{(r_1, \dotsc, r_t)}(S) \defined \set[\big]{A \subseteq S \st \text{$A^{r_i}$ is $L_i$-free for all $1 \leq i \leq t$}},
\end{equation*}
and the maximum-density is defined as
\begin{equation*}
    d_{L_1, \dotsc, L_t}^{(r_1, \dotsc, r_t)}(S) \defined \max\set[\big]{\card{A}/\card{S} \st A \in \cF_{L_1, \dotsc, L_t}^{(r_1, \dotsc, r_t)}(S)}.
\end{equation*}
Write $L(x_1, \dotsc, x_s) = a_1 x_1 + \dotsb + a_s x_s$ and the condition that $A^r$ is $L$-free is equivalent as saying that $A$ has no solution to the \emph{diagonal} equation
\begin{equation*}
    a_1 x_1^r + \dotsb + a_s x_s^r = 0.
\end{equation*}
Or next result says that these constraints are asymptotically independent for different values of $r$.
Note that we are also allowing $r$ to be negative.

\begin{theorem}
\label{thm:general-independence-fq}
For every $L_1, \dotsc, L_t \in \Lin_{\geq 3}$, and distinct non-zero integers $r_1, \dotsc, r_t$, we have that
\begin{equation}
\label{eq:general-independence-fq}
     d_{L_1, \dotsc, L_t}^{(r_1, \dotsc, r_t)}(\FF_q) = d_{L_1}(\FF_q) \dotsb d_{L_t}(\FF_q) + o(1) = q^{-1} \log_2 \card{\cF_{L_1, \dotsc, L_t}^{(r_1, \dotsc, r_t)}(\FF_q)}
\end{equation}
when $q \to \infty$, as long as $L_1, \dotsc, L_t$ are non-degenerate in $\FF_q$ and $r_1, \dotsc, r_t$ are non-zero in $\FF_q$.
\end{theorem}

We prove \Cref{thm:general-independence-fq} via the container\footnote{While we could have used the hypergraph containers results of Saxton and Thomason~\cite{Saxton2015-cr} and of Balogh, Morris, and Samotij~\cite{Balogh2015-al}, our family of containers will be built in explicitly form.} formulation \Cref{thm:general-container} in \Cref{sec:indep}.
Using our version of the wrapper theorem \Cref{thm:polywrapper}, we show that a set that has few solutions to a linear form $L_i \in \Lin(\FF_q)_{\geq 3}$ can be well approximated by generalised Bohr sets known as \emph{wrappers}.
Together with Green's arithmetical removal lemma \Cref{lem:green-removal}, we show that the density of these wrappers cannot be much larger than $d_{L_i}(\FF_q)$.
However, wrappers have their Fourier coefficients controlled in the sense that the Wiener norm of their indicator function is small.
This allows us, using variant of the Weil bound for a generalised character sum in \Cref{cor:cochrane-pinner}, to show that the density of the intersection of these wrappers\footnote{We note that the use of character sums together with wrappers to control intersection sizes is also a key feature in the works of Semchankau~\cite{Semchankau2025-nt}, Semchankau and Shkredov~\cite{Semchankau2026-dq}, Chow, Lim and Mudgal~\cite{Chow2025-hs}, and Shao~\cite{Shao2026-ti}.} is close to $d_{L_1}(\FF_q) \dotsb d_{L_t}(\FF_q)$.
For \Cref{thm:density-independence-fq}, the Weil bound (\Cref{prop:weil}) for Kloosterman  sums suffices.
Since there are not to many wrappers, this suits us with a family of containers which gives the enumeration in \eqref{eq:general-independence-fq}.

The lower bound for \Cref{thm:general-independence-fq}, and consequently of \Cref{thm:density-equation-inverse-fq,thm:density-independence-fq,thm:counting-independence-fq,thm:counting-equation-inverse-fq} is provided in \Cref{sec:constructions}.
The main technical difficulty is to establish that there are $L$-free subsets of $\FF_q$ which have density close to optimal and whose indicator function has suitably small Wiener norm.


\subsection{Product-free sets with linear constraints}

A set $A$ is \emph{product-free} if there are no $x,y,z \in A$ with $x y = z$.
How large can a set $A$ be if it is \emph{simultaneously} sum-free and product-free?
This question was considered first in the work of Kurlberg, Lagarias, and Pomerance~\cite{Kurlberg2012-wo}.
For a subset $S \subseteq R$ of a ring, consider the family
\begin{align*}
    \cF_{\Pi}(S) &\defined \set[\big]{A \subseteq S \st \text{$A$ is product-free}}
\end{align*}
and also $\cF_{\Sigma,\Pi}(S) \defined \cF_\Sigma(S) \cap \cF_\Pi(S)$ for the family of subsets of $S$ that are simultaneously sum-free and product-free.
As usual, define the densities
\begin{align*}
    d_\Pi(S) &\defined \max \set[\big]{\card{A}/\card{S} \st \text{$A \subseteq S$ is product-free } }, \\
    d_{\Sigma,\Pi}(S) &\defined \max \set[\big]{\card{A}/\card{S} \st  \text{$A \subseteq S$ is sum-free and product-free }}.
\end{align*}
Kurlberg, Lagarias, and Pomerance~\cite{Kurlberg2012-wo,Kurlberg2013-bk} studied these quantities in $\ZZ/N\ZZ$.
Define the function $f(x) \defined (\loglog x)^{1-(e/2)\log 2} (\logloglog x)^{1/2}$.
They proved that
\begin{equation*}
    d_{\Pi}(\ZZ/N\ZZ) \leq 1 - \frac{\kappa}{f(N)}, \qquad d_{\Sigma,\Pi}(\ZZ/N\ZZ) \leq \frac{1}{2} - \frac{\kappa'}{f(N)},
\end{equation*}
for some $\kappa, \kappa' > 0$ and for large enough $N$.
At the same time, they have shown that the order of the function $f$ is optimal in this theorem.
That is, there are $\kappa'', \kappa''' > 0$ for which each of the following inequalities hold for infinitely many $N$:
\begin{equation*}
    d_{\Pi}(\ZZ/N\ZZ) \geq 1 - \frac{\kappa''}{f(N)}, \qquad d_{\Sigma,\Pi}(\ZZ/N\ZZ) \geq \frac{1}{2} - \frac{\kappa'''}{f(N)}.
\end{equation*}
From the characterisation of $d_\Sigma(G)$ for finite commutative groups, the value of $d_\Sigma(\ZZ/N\ZZ)$ is known exactly for all $N$ and can be read of the prime factorisation of $N$.
Indeed, we have
\begin{equation*}
    d_\Sigma(\ZZ/N\ZZ) = \begin{cases}
    1/3 - 1/3N & \text{if every prime $p \divides N$ satisfy $p \eqmod 1 \mod{3}$}, \\
    1/3 + 1/3p & \text{if $p$ is the least prime $p \divides N$ with $p \eqmod 2 \mod{3}$}, \\
    1/3 & \text{if $3 \divides N$ and no prime $p \divides N$ satisfies $p \eqmod 2 \mod{3}$}.
    \end{cases}
\end{equation*}
Note that $d_\Sigma(\ZZ/N\ZZ) \leq 1/2$ and equality occur for even $N$, as the set of odds is sum-free.
The results above suggest that $d_{\Sigma,\Pi}(\ZZ/N\ZZ) = d_\Sigma(\ZZ/N\ZZ)d_\Pi(\ZZ/N\ZZ) + o(1)$ holds at least along some subsequences.
We now show that this asymptotic independence of sum-freeness and product-freeness occurs over all finite fields.

\begin{theorem}
\label{thm:sum-product}
We have that $d_{\Sigma,\Pi}(\FF_q) = d_\Sigma(\FF_q)d_\Pi(\FF_q) + o(1)$ when $q \to \infty$.
\end{theorem}

Note that $d_\Pi(\FF_q) = (1 - 1/q) d_\Pi(\FF_q^\times)$ and that $\cF_\Pi(\FF_q)$ is in bijection to $\cF_\Sigma(\ZZ/(q-1)\ZZ)$, hence $d_\Pi(\FF_q^\times) = d_\Sigma(\ZZ/(q-1)\ZZ)$.
In principle, $d_{\Sigma,\Pi}(\FF_q)$ may not converge except along some subsequences since $d_\Sigma(\ZZ/(q-1)\ZZ)$ may depends wildly on the prime factorisation of $q - 1$.
However, the dependency is well behaved.
If $q$ is odd then $q - 1$ is even, and hence $d_\Pi(\FF_q) = 1/2$. 
This gives that $d_{\Sigma,\Pi}(\FF_q) = d_{\Sigma}(\FF_q)/2 + o(1)$ as $q \to \infty$.
For instance, this gives $d_{\Sigma,\Pi}(\FF_{7^d}) = 1/7 + o(1)$.

If $A \subseteq \FF_p$ is sum-free, then $A + A$ is disjoint from $A$, and if it is product-free, then $A \cdot A$ is disjoint from $A$.
If both conditions are simultaneously satisfied, one may estimate $d_{\Sigma,\Pi}(\FF_p)$ via a \emph{sum-product} theorem as follows.
Suppose that for every set $A \subseteq \FF_p$ with $\card{A} = \alpha p$ one has $\max\set{\card{A + A}, \card{A \cdot A}} \geq f(\alpha) p$ for some function $f(\alpha)$, then $\beta = d_{\Sigma,\Pi}(\FF_p)$ satisfies $\beta + f(\beta) \leq 1$.
Shao~\Cite{Shao2026-ti} has recently determined exactly the optimal function $f$ for this sum-product theorem in dense subsets of $\FF_p$.
His result gives $\max\set{\card{A + A}, \card{A \cdot A}} \geq (4\alpha - o(1)) p$ if $1/8 \leq \alpha \leq 1/4$, which leads to the bound $d_{\Sigma,\Pi}(\FF_p) \leq 1/5 + o(1)$.
We now establish via \Cref{thm:sum-product} that
\begin{equation*}
    d_{\Sigma,\Pi}(\FF_p) = 1/6 + o(1),
\end{equation*}
which indicates that additional information beyond the sizes of $A + A$ and $A \cdot A$ must be used to obtain the optimal density of sum-free product-free sets.

Now turn to the case $q = 2^d$, so $d_\Pi(\FF_q)$ depends on the prime factorisation of $2^d - 1$.
We have indeed that $d_\Pi(\FF_q) = 1/3 + o(1)$ unless $2^d - 1$ has a prime factor $p$ with $p \eqmod 2 \mod{3}$.
We show in \Cref{prop:produc-liminfsup} that
\begin{equation}
\label{eq:liminfsup}
    \frac{1}{3} = \liminf_{d \to \infty} d_\Pi(\FF_{2^d}) \leq \limsup_{d \to \infty} d_\Pi(\FF_{2^d}) = \frac{2}{5}. 
\end{equation}
Hence, $d_{\Sigma,\Pi}(\FF_{2^d})$ is asymptotically bounded between $1/6$ and $1/5$ and these values are attained in the limit for certain subsequences.

While we have shown that sum-freeness and product-freeness are asymptotically independent in, much more is true. 
Consider the following family
\begin{equation*}
    \cF_{L_1, \dotsc, L_t,\Pi}^{(r_1, \dotsc, r_t)}(S) \defined \set[\big]{A \subseteq S \st \text{$A^{r_i}$ is $L_i$-free for all $1 \leq i \leq t$, $A$ is product-free}},
\end{equation*}
with the corresponding density function
\begin{equation*}
    d_{L_1, \dotsc, L_t,\Pi}^{(r_1, \dotsc, r_t)}(S) \defined \max\set[\big]{\card{A}/\card{S} \st A \in \cF_{L_1, \dotsc, L_t,\Pi}^{(r_1, \dotsc, r_t)}(S)}.
\end{equation*}
Our final result is a generalisation of \Cref{thm:sum-product}.

\begin{theorem}
\label{thm:general-independence-fq-product}
For every $L_1, \dotsc, L_t \in \Lin_{\geq 3}$, and distinct non-zero integers $r_1, \dotsc, r_t$, we have that
\begin{equation}
\label{eq:general-independence-fq-product}
     d_{L_1, \dotsc, L_t, \Pi}^{(r_1, \dotsc, r_t)}(\FF_q) = d_{L_1}(\FF_q) \dotsb d_{L_t}(\FF_q)d_\Pi(\FF_q^\times) + o(1) = q^{-1} \log_2 \card{\cF_{L_1, \dotsc, L_t, \Pi}^{(r_1, \dotsc, r_t)}(\FF_q)}
\end{equation}
when $q \to \infty$, as long as $L_1, \dotsc, L_t$ are non-degenerate in $\FF_q$ and $r_1, \dotsc, r_t$ are non-zero in $\FF_q$.
\end{theorem}

The novelty in the proof of \Cref{thm:general-independence-fq-product} compared with \Cref{thm:general-independence-fq} is the full use of the \emph{mixed character sums} in \Cref{cor:cochrane-pinner}.
Indeed, we simultaneously consider the Fourier decomposition of a subset $A \subseteq \FF_q$ into additive characters of $\FF_q$ and multiplicative characters of $\FF_q^\times$.
We then apply the wrapper machinery to each group separately and combine the results the character sum.

\begin{remark}
The density of $L$-free sets is the same as the density of sets with no solutions to $L(x) = b$, given that $L \in \Lin$ is not translation invariant.
Indeed, $A$ is $L$-free if and only if the translation $A + b/(a_1 + \dotsb + a_s)$ has no solution to $L(x) = a_1 x_1 + \dotsb + a_s x_s = b$.
Therefore, every asymptotic result above, such as \Cref{thm:general-independence-fq} or \Cref{thm:general-independence-fq-product}, still holds if the condition that $A^{r_i}$ is $L_i$-free gets replaced by the condition that $A^{r_i}$ has no solution to $L_i(x) = b_i$ for some $b_1, \dotsc, b_t \in \FF_q$, given the same non-degeneracy conditions. 
\end{remark}


\section*{Methodology}

This work originated during a Supervised Research course in Spring 2026 at Cornell University pursued by the first author and supervised by the second author.
All the mathematical results were obtained by the authors, who also wrote the text in full.
Language models were used for literature search and for finding minor mistakes, typos, and inconsistencies.


\section{Preliminaries}


\subsection{Asymptotic notation}

We use Vinogradov notation, namely, we write $X \ll Y$ to mean that there is some absolute constant $C > 0$ such that $\abs{X} \leq C Y$.
We write $X \ll_z Y$ when the leading constant $C$ above is a function of $z$.
More generally, we may write $X \ll_{z_1, \dotsc, z_t} Y$ if $\abs{X} \leq C Y$ for some $C = C(z_1, \dotsc, z_t) > 0$.
We also write $X = O_z(Y)$ for $X \ll_z Y$.
We write $f(x) = o(g(x))$ as $x \to \infty$ when the the limit  $\lim_{x \to \infty}f(x)/g(x)$ exists and is equal to $0$.


\subsection{Fourier analysis}

Let $G$ be a a finite commutative group.
Instead of working with subsets $A \subseteq G$, we work with their indicator functions $\ind_A \from G \to \set{0,1}$ instead, seen as elements of the space of functions $\rL(G) \defined \set{f \from G \to \CC}$.
Formally, we equip the $G$ with the uniform probability measure $S \mapsto \card{S}/\card{G}$ for $S \subseteq G$.
From this measure in $G$, we consider the space of functions $\rL(G)$ with the $L^r$-norms 
\begin{equation*}
    \norm{f}_r \defined \p[\Big]{\frac{1}{\card{G}} \sum_{x \in G} \abs{f(x)}^r}^{1/r},
\end{equation*}
for all $r \geq 1$, with the $L^\infty$-norm
\begin{equation*}
    \norm{f}_\infty \defined \max\set[\big]{\abs{f(x)} \st x \in G},
\end{equation*}
and with inner product defined for all $f, g \in \rL(G)$ as
\begin{equation*}
    \inner{f}{g} \defined \frac{1}{\card{G}} \sum_{x \in G}f(x)\conj{g(x)}.
\end{equation*}
Let $\SS \subseteq \CC$ be the unit complex circle.
A \emph{character} in $G$ is a group homomorphism $\chi \from G \to (\SS, \cdot)$.
Let $\fourier{G} = \set{\chi \st \text{$\chi$ is a character}}$ be the set of characters of $G$, and recall that $\fourier{G}$ is a group with pointwise multiplication and identity element being the \emph{trivial character} $\chi_0$ satisfying $\chi_0(x) = 1$ for all $x \in G$.
Note that $\fourier{G} \isomorphic G$ as an commutative group.
Given a function $f \in \rL(G)$, its Fourier transform $\fourier{f} \from \fourier{G} \to \CC$ is given by
\begin{equation*}
    \fourier{f}(\chi)
    \defined \inner{f}{\chi} = \frac{1}{\card{G}}\sum_{x \in G} f(x)\conj{\chi(x)}.
\end{equation*}
Note that every function $\eta \from \fourier{G} \to \CC$ is the Fourier transform of a function $G \to \CC$, namely, the function $x \mapsto \sum_{\chi \in \fourier{G}} \eta(\chi) \chi$.
Thus, we have $\rL(\fourier{G}) = \set{\fourier{f} \st f \in \rL(G)}$.
It is convenient to consider $\fourier{G}$ with the counting measure $S \mapsto \card{S}$ for $S \subseteq \fourier{G}$.
Thus, we equip $\rL(\fourier{G})$ with the following norms and inner product
\begin{align*}
    \norm{\fourier{f}}_r \defined \p[\Big]{\sum_{\chi \in \fourier{G}} \abs[\big]{\fourier{f}(\chi)}^r}^{1/r}, \qquad 
    \norm{\fourier{f}}_\infty \defined \max\set[\big]{\abs[\big]{\fourier{f}(\chi)} \st \chi \in \fourier{G}},\qquad
    \inner{\fourier{f}}{\fourier{g}} \defined \sum_{\chi \in \fourier{G}}\fourier{f}(\chi)\conj{\fourier{g}(\chi)}.
\end{align*}
for $r \geq 1$ and $f, g \from G \to \CC$.

Following the same measure convention, the convolution in $\rL(G)$ is defined as,
\begin{align*}
    (f \conv g)(z) &\defined \frac{1}{\card{G}} \sum_{x \in G} f(x) g(z-x) = \frac{1}{\card{G}}\sum_{\substack{x + y = z \\ x,y \in G}} f(x)g(y).
\end{align*}
while the convolution in $\rL(\fourier{G})$ is
\begin{align*}
    (\fourier{f} \conv \fourier{g})(\chi) &\defined \sum_{\chi \in \fourier{G}} \fourier{f}(\eta) g(\chi \eta^{-1}) = \sum_{\substack{\eta \psi = \chi \\ \eta, \psi \in \fourier{G}}} \fourier{f}(\eta)\fourier{g}(\psi).
\end{align*}
Our choice or normalisation has the following useful consequences.
For every $f, g \from G \to \CC$, $x \in G$, $\chi \in \fourier{G}$, $A \subseteq G$ with density $\alpha \defined \card{A}/\card{G}$, we have the \emph{Fourier decomposition}
\begin{equation*}
    f = \sum_{\chi \in \fourier{G}} \fourier{f}(\chi)\, \chi.
\end{equation*}
Furthermore, we have the \emph{convolution identities} $\fourier{f \conv g} = \fourier{f} \cdot \fourier{g}$ and $\fourier{f} \conv \fourier{g} = \fourier{f \cdot g}$, the \emph{Parseval identity} $\inner{f}{g} = \inner{\fourier{f}}{\fourier{g}}$, and the \emph{Planceherel identity} $\norm{f}_2 = \norm{\fourier{f}}_2$ for every $f,g \in \rL(G)$.
We also have the \emph{Young inequalities} $\norm{f \conv g}_1 \leq \norm{f}_1 \norm{g}_1$ and $\norm{\fourier{f} \conv \fourier{g}}_1 \leq \norm{\fourier{f}}_1 \norm{\fourier{g}}_1$.

Since the convolutions are associative, repeated convolution in $\rL(G)$ satisfies
\begin{align*}
    (f_1 \conv \dotsb \conv f_k)(z) &= \frac{1}{\card{G}^{k-1}} \sum_{x_1, \dotsc, x_{k-1} \in G}f_1(x_1) \dotsb f_{k-1}(x_{k-1})f_k(z - x_1 - \dotsb - x_{k-1}) \\
    &= \frac{1}{\card{G}^{k-1}} \sum_{\substack{x_1 + \dotsb + x_k = z \\ x_1,\dotsc,x_k \in G}}f_1(x_1) \dotsb f_k(x_k),
\end{align*}
while by Fourier decomposition followed by repeated convolution identity, we have
\begin{equation*}
    (f_1 \conv \dotsb \conv f_k) (z) = \sum_{\chi \in \fourier{G}} (\widehat{f_1 \conv \dotsb \conv f_k})(\chi)\,\chi(z) = \sum_{\chi \in \fourier{G}} \fourier{f_1}(\chi) \dotsb \fourier{f_k}(\chi)\,\chi(z).
\end{equation*}
Also note that the triangle inequality gives
\begin{equation*}
    \abs{\fourier{f}(\chi)} = \frac{1}{\card{G}}\abs[\Big]{\sum_{x \in G}f(x)\conj{\chi(x)}} \leq \frac{1}{\card{G}} \sum_{x \in G}\abs{f(x)} = \norm{f}_1,
\end{equation*}
for all $f \in \rL(G)$ and $\chi \in \fourier{G}$.
The indicator of a set $A \subseteq G$ satisfies the identities:
\begin{equation*}
    \fourier{\ind}_A(\chi_0) = \norm{\ind_A}_2^2 = \norm{\fourier{\ind}_A}_2^2 = \alpha.
\end{equation*}


\subsection{Finite fields}

Let $\FF_q$ be a finite field of characteristic $p$ and $q = p^d$.
$\FF_q$ has two commutative group structures, the \emph{additive} $(\FF_q, +)$ and the \emph{multiplicative} $(\FF_q^\times, \times)$.

We write $\Tr = \Tr_{\FF_q / \FF_p} \from \FF_q \to \FF_p$ for the trace map, and we fix a standard non-trivial additive character $\psi \from (\FF_q, +) \to (\CC^{\times}, \cdot)$, defined as
\begin{equation*}
    \psi(x) \defined \exp\p[\big]{2 \pi i \Tr(x)/p}.
\end{equation*}
Parametrise the group $\fourier{\FF}_q$ of additive characters as $\fourier{\FF}_q = \set{\psi_a \from a \in \FF_q}$, where for $a \in \FF_q$,
\begin{equation*}
    \psi_a(x) \defined \psi(ax).
\end{equation*}
We then have $\psi_a \psi_b = \psi_{a+b}$ and that $\fourier{\FF}_q \isomorphic \FF_q$.

The multiplicative group $(\FF_q^\times, \times)$ is cyclic with $q - 1$ elements.
Fix a \emph{primitive root} $\gamma$ in $\FF_q^\times$ define a standard non-trivial multiplicative character $\chi \from (\FF_q^\times, \times) \to (\CC^\times, \cdot)$, defined as
\begin{equation*}
    \chi(\gamma^k) \defined \exp\p[\big]{2 \pi i k / (q-1)}.
\end{equation*}
We parametrise the group $\fourier{\FF_q^\times} \defined \set{\chi_a \st a \in \FF_q^\times}$, where fore $a = \gamma^m$, we define
\begin{equation*}
    \chi_a(x) \defined \chi(x^m).
\end{equation*}
In other words, $\chi_{\gamma^m}(x) \defined \chi(x^m)$, which gives that $\chi_a \chi_b = \chi_{a b}$, the trivial character being $\chi_1$.


\subsection{The Wiener norm}

The \emph{Wiener norm} of a function $f \from G \to \CC$ is the $L^1$-norm of its Fourier transform, namely
\begin{equation*}
    \norm{f}_\wiener \defined \norm{\fourier{f}}_{1} = \sum_{\chi \in \fourier{G}} \abs{\fourier{f}(\chi)}.
\end{equation*}
We write $\norm{f}_{\wiener(G)}$ when we want to emphasise that we are taking $f$ as a function in $\rL(G)$.
The Wiener norm also goes by the name of \emph{algebra norm}, since for every $f, g \in \rL(G)$, we have
\begin{equation}
\label{eq:algebra-norm}
    \norm{f g}_\wiener = \norm{\fourier{f g}}_1 = \norm{\fourier{f} \conv \fourier{g}}_1 \leq \norm{\fourier{f}}_1 \norm{\fourier{g}}_1 = \norm{f}_\wiener \norm{g}_\wiener
\end{equation}
by the Young inequality.
The following inequality for the product of $k \geq 2$ functions is also used
\begin{align}
\label{eq:wiener-convolution}
\begin{split}
    \norm{f_1 \conv \dotsb \conv f_k}_\wiener 
    &= \sum_{\chi \in \fourier{G}} \abs[\big]{(\fourier{f_1 \conv \dotsb \conv f_k})(\chi)}
    = \sum_{\chi \in \fourier{G}} \abs{\fourier{f_1}(\chi)} \dotsb \abs{\fourier{f_k}(\chi)} \\
    &\leq \norm{f_3}_1 \dotsb \norm{f_k}_1 \sum_{\chi \in \fourier{G}} \abs{\fourier{f_1}(\chi)}\abs{\fourier{f_2}(\chi)} \leq \norm{f_1}_2 \norm{f_2}_2 \norm{f_3}_1 \dotsb \norm{f_k}_1,
\end{split}
\end{align}
where the last step we applied the Cauchy-Schwarz inequality.
If we take $f_1 = \ind_A$ and $f_2 = \card{G} \ind_{\set{0}}$, we have $f_1 \conv f_2 = \ind_A$ and \eqref{eq:wiener-convolution} gives the bound
\begin{equation}
\label{eq:wiener-cardinality}
    \norm{\ind_A}_\wiener \leq \norm{f_1}_2 \norm{f_2}_2 = \card{A}^{1/2}.
\end{equation}

For an automorphism $\alpha \from G \to G$ we have $\fourier{(f \compo \alpha)}(\chi) = \fourier{f}(\chi \compo \alpha^{-1})$, thus $\norm{f \compo \alpha}_\wiener = \norm{f}_\wiener$.
Thus translation, negation, and dilation by an $k \in \ZZ$ with $\gcd(k, \card{G}) = 1$ preserve the Wiener norm.

Let $G$ and $H$ be finite commutative groups and $f \in \rL(G)$ and $g \in \rL(H)$.
Then the function $f \otimes g \in \rL(G \times H)$ defined by $(f \otimes g)(x,y) = f(x) g(y)$ satisfies
\begin{equation}
\label{eq:wiener-product}
    \norm{f \otimes g}_{\wiener(G \times H)} = \norm{f}_{\wiener(G)} \norm{g}_{\wiener(H)}.
\end{equation}
Indeed, this follows from the identity $\fourier{(f \otimes g)}(\chi \otimes \eta) = \fourier{f}(\chi) \fourier{g}(\eta)$ for every $\chi \in \fourier{G}$, $\eta \in \fourier{H}$.

\begin{lemma}
\label{lem:wiener-restriction}
Let $H$ be a subgroup of a finite commutative group $G$.
For every $f \from G \to \CC$ we have
\begin{equation*}
    \norm{f \!\restriction_H }_{\wiener(H)} \leq \norm{f}_{\wiener(G)},
\end{equation*}
where $f \!\restriction_H \from H \to \CC$ is the restriction of $f$ to the subgroup $H$.
\end{lemma}
\begin{proof}
For every character $\chi \in \fourier{G}$, we have that $\chi \!\restriction_H$ is a character in $H$.
We start from the Fourier decomposition $f = \sum_{\chi \in \fourier{G}} \fourier{f}(\chi)\, \chi$ and group the characters in $\fourier{G}$ by their restriction into $\fourier{H}$ as follows:
\begin{equation*}
    f\!\restriction_H = \sum_{\chi \in \fourier{G}} \fourier{f}(\chi)\, \chi\!\restriction_H = \sum_{\eta \in \fourier{H}} \sum_{\substack{\chi \in \fourier{G} \\ \chi \restriction_H = \eta}} \fourier{f}(\chi)\, \eta.
\end{equation*}
Therefore, the Fourier coefficients of the restriction are precisely given by
\begin{equation*}
    \fourier{f\!\restriction_H}(\eta) =  \sum_{\substack{\chi \in \fourier{G} \\ \chi \restriction_H = \eta}} \fourier{f}(\chi).
\end{equation*}
Finally, the triangle inequality gives
\begin{equation*}
    \norm{f\!\restriction_H}_{\wiener(H)} = \sum_{\eta \in \fourier{H}}\abs{\fourier{f\!\restriction_H}(\eta)}
    \leq \sum_{\eta \in \fourier{H}} \sum_{\substack{\chi \in \fourier{G} \\ \chi \restriction_H = \eta}} \abs{\fourier{f}(\chi)} = \norm{f}_{\wiener(G)}. \qedhere
\end{equation*}
\end{proof}

\begin{lemma}
\label{lem:wiener-surjective}
Let $\pi \from G \to H$ be a surjective homomorphism of finite commutative groups.
Then for every $f \from H \to \CC$, we have
\begin{equation*}
    \norm{f \compo \pi}_{\wiener(G)} = \norm{f}_{\wiener(H)}.
\end{equation*}
\end{lemma}
\begin{proof}
Note that for each character $\chi \in \fourier{H}$, the map $\chi \compo \pi \from G \to \CC$ is also a character.
Moreover, as $\pi$ is surjective, all the characters of the form $\chi \compo \pi$, for $\chi \in \fourier{H}$, are distinct.
From the Fourier decomposition $f = \sum_{\chi \in \fourier{H}} \fourier{f}(\chi)\, \chi$ we get
\begin{equation*}
    f \compo \pi = \sum_{\chi \in \fourier{H}} \fourier{f}(\chi) \,\chi \compo \pi.
\end{equation*}
Hence $\fourier{(f \compo \pi)}(\chi \compo \pi) = \fourier{f}(\chi)$ and $\fourier{(f \compo \pi)}(\eta) = 0$ if $\eta$ is not of the form $\chi \compo \pi$ for some $\chi \in \fourier{H}$.
Thus
\begin{equation*}
    \norm{f \compo \pi}_{\wiener(G)} = \sum_{\eta \in \fourier{G}} \abs{\fourier{(f \compo \pi)}(\eta)} = \sum_{\chi \in \fourier{H}} \abs{\fourier{f}(\chi)} = \norm{f}_{\wiener(H)}.
    \qedhere
\end{equation*}
\end{proof}

The following quantity will be relevant for us
\begin{equation*}
    \cB(G) \defined \sup\set[\big]{ \norm{\ind_{\chi^{-1} (A)}}_\wiener \st \text{$\chi \in \fourier{G}$ and $A \subseteq \SS$ is an arc}}.
\end{equation*}
Indeed, the parameter $\cB(G)$ controls the the Wiener norm of the indicator of the inhomogeneous Bohr set $\set{x \in G \st \chi(x) \in A}$.
It turns out that $\cB(G)$ is closely related to the parameter
\begin{equation*}
    \cI(G) \defined \max\set[\big]{\norm{\ind_P}_\wiener \st \text{$P \subseteq G$ is an arithmetic progression}},
\end{equation*}
where a set is an arithmetic progression if it is of the form $\set{a k + b \st 1 \leq k \leq L}$ for some $a, b \in G$ and integer $L$.
Since the Wiener norm is invariant under translation, it suffices to consider $b = 0$.
Furthermore, from \Cref{lem:wiener-restriction} we get that $\cI(H) \leq \cI(G)$ when $H$ is a subgroup of $G$.
Indeed, any arithmetic progression $P \subseteq H$ is also an arithmetic progression in $G$, thus $\norm{\ind_P}_{\wiener(H)} \leq \norm{\ind_P}_{\wiener(G)} \leq \cI(G)$.
As a consequence, we have that
\begin{equation*}
    \cI(\ZZ/N\ZZ) = \max\set[\big]{\norm{\ind_{\set{0,1,\dotsc,L}}}_\wiener \st 0 \leq L \leq N}.
\end{equation*}

We will now show that this parameter is controlled by the \emph{exponent} $\exp(G)$ of the group $G$, namely, the smallest $N \in \ZZ_{\geq 1}$ such that $Nx = 0$ for all $x \in G$.
Since $G$ is finite and commutative, $N$ is also the size of the largest cyclic subgroup of $G$.

\begin{lemma}
If $G$ is finite commutative group with exponent $N = \exp(G)$, then $\cB(G) = \cI( \ZZ / N \ZZ)$.
\end{lemma}
\begin{proof}
Let $\chi \in \fourier{G}$ be an element of order $d$ in $\fourier{G}$.
Then $d$ divides $N$ and $\chi(G) \subseteq \CC$ is the group of $d$-th roots of unity $\SS_d$, isomorphic to $\ZZ/d\ZZ$.
Thus $\chi \from G \to \SS_d$ is a surjective group homomorphism.
For any arc $A \subseteq \SS$, the subset $P = A \cap \SS_d$ is an arithmetic progression in $\SS_d$, so $\norm{\ind_P}_\wiener \leq \cI(\ZZ/d\ZZ)$.
Note that $\ind_{\chi^{-1}(A)} = \ind_P \compo \chi$, so \Cref{lem:wiener-surjective} gives
\begin{equation*}
    \norm{\ind_{\chi^{-1}(A)}}_{\wiener(G)} = \norm{\ind_P \compo \chi}_{\wiener(G)} = \norm{\ind_P}_{\wiener(\SS_d)} \leq \cI(\ZZ/d\ZZ) \leq \cI(\ZZ/N\ZZ),
\end{equation*}
the last inequality being due to \Cref{lem:wiener-restriction} since $\ZZ/d\ZZ$ is a subgroup of $\ZZ/N\ZZ$.
This gives the upper bound $\cB(G) \leq \cI(\ZZ/N\ZZ)$.

For the lower bound, note that $G$ has a subgroup $H$ isomorphic to $\ZZ/N\ZZ$ since $G$ has an element of order $N$.
Let $\pi \from G \to H$ be a surjective homomorphism.
Let $I \subseteq H \isomorphic \ZZ/N\ZZ$ be an interval of the form $\set{0,1,\dotsc,L}$.
Take the character $\chi(x) = \exp(2 \pi i x/N)$ in $\fourier{H}$ and note that $\chi \compo \pi$ is a character in $\fourier{G}$.
Now the circular arc $A = \set{\exp(2 \pi i \theta/N) \st 0 \leq \theta \leq L} \subseteq \SS$, and note that $\chi^{-1}(A) = I$ in $H$ and that $(\chi \compo \pi)^{-1}(A) \cap H = I$.
Therefore, from \Cref{lem:wiener-restriction}, we have
\begin{equation*}
    \cB(G) \geq \norm[\big]{\ind_{(\chi \compo \pi)^{-1}(A)}}_{\wiener(G)}
    \geq \norm[\big]{\ind_{(\chi \compo \pi)^{-1}(A)} \restriction_H}_{\wiener(H)}
    = \norm[\big]{\ind_{I}}_{\wiener(H)}.
\end{equation*}
Since this holds for every circular interval $I$, we have $\cB(G) \geq \cI(H) = \cI(\ZZ/N\ZZ)$, as claimed.
\end{proof}

The following proposition is a standard bound on $\cI(\ZZ/n\ZZ)$.

\begin{proposition}
\label{prop:wienerinterval}
For every $n \geq 2$, we have $\cI(\ZZ/n\ZZ) \leq 2\log n$.
\end{proposition}
\begin{proof}
Since the Wiener norm is invariant under translations and dilations, it is enough to consider the Wiener norm of $\ind_A$ where $A = \set{1,\dots,L} \subseteq \ZZ/n\ZZ$ is an interval.
If $n = 2$, then it can be seen that $\norm{\ind_A}_\wiener \leq 1$ as either $A = \ZZ/n\ZZ$ or $L \leq 1$.
Parametrise $\fourier{\ZZ/n\ZZ} = \set{\psi_a \st -(n-1)/2 \leq a \leq (n-1)/2}$.
Write $\beta \defined \exp(-2\pi i a /n)$ and observe that $\fourier{\ind}_A(\psi_a)$ is the following geometric sum
\begin{align*}
    \fourier{\ind}_A(\psi_a) &= \frac{1}{n} \sum_{x \in A} \psi_a(-x) = \frac{1}{n}\sum_{1 \leq k \leq L} \beta^{k} = \frac{1}{n} \cdot \frac{\beta^{L+1} - \beta}{\beta - 1} = \frac{\beta^{(L+1)/2}}{n} \cdot \frac{\sin(\pi a (L+1)/n)}{\sin(\pi a / n)}.
\end{align*}
Now, note that $\abs{\sin(\pi x)} \geq 2\abs{x}$ for $0 \leq \abs{x} \leq 1/2$, and thus
\begin{equation*}
    \abs{\fourier{\ind}_A(\psi_a)} \leq \frac{1}{n \abs{\sin(\pi a / n)}} \leq \frac{1}{2n \abs{a/p}} = \frac{1}{2\abs{a}}.
\end{equation*}
Therefore, for $n \geq 3$ we have
\begin{equation*}
    \norm{\ind_A}_{\wiener} = \sum_{a \in \ZZ/n\ZZ} \abs{\fourier{\ind}_A(\psi_a)} \leq \frac{L}{n} + 2 \sum_{a=1}^{(n-1)/2} \frac{1}{2a} \leq 1 + 1 + \log((n-1)/2) \leq 2\log n,
\end{equation*}
where we used the fact that $1 + 1/2 + \dotsc + 1/k \leq 1 + \log k$.
\end{proof}

Only the constant in the bound $\cI(\ZZ/n\ZZ) \leq 2 \log n$ can be improved, but the bound $\cI(\ZZ/n\ZZ) \ll \log n$ is already stronger than what we need.

Since every finite commutative group is a product of cyclic groups and we have
\begin{equation*}
    \exp\p[\big]{(\ZZ/n_1 \ZZ) \times \dotsb \times (\ZZ/n_k \ZZ)} = \lcm(n_1, \dotsc, n_k),
\end{equation*}
we obtain the following corollary.

\begin{corollary}
\label{cor:phi-G}
For every finite commutative group $G$, we have $\cB(G) \leq 2 \log \card{G}$.
Moreover, if $q = p^d$ for a prime $p$ and $n \geq 1$, we have $\cB(\FF_q^n) \leq 2 \log p$.
\end{corollary}


\subsection{Exponential sums}

For $a, b \in \FF_q$, the \emph{Kloosterman sums} are 
\begin{equation*}
    \rK(a,b) \defined \sum_{x \in \FF_q^\times}\psi(ax + bx^{-1}) = \sum_{x \in \FF_q^\times} \psi_a(x)\psi_b(x^{-1}).
\end{equation*}
Note that $\rK(0,0) = q - 1$ and that for $a,b \neq 0$, we have $\rK(a,0) = \rK(0,b) = -1$.
We will make use of the following classical bound of Weil~\cite{Weil1948-gy}.

\begin{proposition}[Weil bound]
\label{prop:weil}
For $a, b \in \FF_q$ with $(a,b) \neq (0,0)$, we have $\abs{\rK(a,b)} \leq 2 q^{1/2}$.
\end{proposition}

The Weil bound implies the following standard mixing estimate.

\begin{proposition}
\label{prop:kloosterman-mixing}
For any sets $A,B\subseteq \FF_q^\times$ with densities $\alpha \defined \card{A}/q$ and $\beta \defined \card{B}/q$, we have
\begin{equation*}
    \abs[\big]{\card{A\cap B^{-1}} - \alpha\beta q} \leq 2q^{1/2} \norm{\ind_A}_\wiener \norm{\ind_B}_\wiener.
\end{equation*}
\end{proposition}
\begin{proof}
By Fourier inversion, we have
\begin{equation*}
    \card{A\cap B^{-1}} = \sum_{x\in \FF_q^\times} \ind_A(x) \ind_B(x^{-1}) = \sum_{a,b\in \FF_q} \fourier{\ind_A}(\psi_a) \fourier{\ind_B}(\psi_b) \rK(a,b).
\end{equation*}
The contribution of the pair $(a,b) = (0,0)$ is exactly $\alpha \beta (q - 1)$. 
Furthermore, the contribution over pairs $(a,b) = (a,0)$ with $a \neq 0$ is
\begin{align*}
    \sum_{a \in \FF_q^\times} \fourier{\ind_{A}}(\psi_a) \fourier{\ind_{B}}(\psi_0) \rK(a,0)
    &= -\beta \sum_{a \in \FF_q^\times} \fourier{\ind_{A}}(\psi_a)
    = -\beta\p[\big]{\ind_A(0) - \beta} = \alpha \beta.
\end{align*}
Hence, by the Weil bounds in \Cref{prop:weil}, we have
\begin{align*}
    \abs[\Big]{\card{A \cap B^{-1}} - \alpha \beta q}
    &\leq \sum_{\substack{a \in \FF_q^\times \\ b \in \FF_q}} \abs[\Big]{\fourier{\ind_{A}}(\psi_a) \fourier{\ind_{B}}(\psi_b) \rK(a, b)} \\
    &\leq 2 q^{1/2} \sum_{ a,b \in \FF_q} \abs[\Big]{\fourier{\ind_{A}}(\psi_a) \fourier{\ind_{B}}(\psi_b)}
    = 2 q^{1/2} \norm{\ind_A}_\wiener \norm{\ind_B}_\wiener. \qedhere
\end{align*}
\end{proof}

While the Kloosterman sums above are sufficient for a mixing of two sets, or more general independence results require a more involved exponential sum.
We use a version of the Weil bound for mixed characters sums of a Laurent polynomial.
The formulation below was adapted from Cochrane and Pinner~\cite{Cochrane2006-ny}, but this bound was already obtained by Perel'muter~\cite{Perel-muter1969-rg} (with a worse constant), and also appear in Katz and Laumon~\cite{Katz1985-fc}, and in Castro and Moreno~\cite{Castro2000-lh}.
Denote by $\FF_q(x)$ the field of rational functions with coefficients in $\FF_q$.

\begin{theorem}
\label{thm:cochrane-pinner}
Let $\FF_q$ be a finite field with characteristic $p$, let $\psi_a$ be an non-trivial additive character and $\chi_b$ be a multiplicative character.
Let $f \in \FF_q[x,x^{-1}]$ be a Laurent polynomial with $f(x) = \sum_{k \in R} a_k x^k$ for $a_k \neq 0$ and some finite $R \subseteq \ZZ$.
If there is no $h \in \FF_q(x)$ and $c \in \FF_q$ with $a f = h^p - h + c$.
Then for $c \defined (\max(\set{0} \cup R) + \max(\set{0} \cup -R))$, we have
\begin{equation*}
    \abs[\Big]{\sum_{x \in \FF_q^\times} \psi_a(f(x)) \chi_b(x) } \leq  c q^{1/2}.
\end{equation*}
\end{theorem}

The reason that $a f = h^p - h + c$ must be avoided to obtain a square-root cancellation is that if this identity holds, then we have $\Tr(af(x)) = \Tr(h^p(x) - h(x) + c) = \Tr(c)$ since $\Tr(x^p) = \Tr(x)$, so $\psi_a(f(x)) = \exp(2\pi \Tr(c)/ p)$ is a constant.
We deal with this constraint in the following lemma.

\begin{lemma}
\label{lem:degenerate}
Let $f(x) \in \FF_q[x,x^{-1}]$ be a non-constant Laurent polynomial $f(x) = \sum_{r \in R} a_r x^r$ with $a_r \neq 0$ and where $r$ is coprime to $p$ for every $r \in R$.
Then there is no $h \in \FF_q(x)$ and no $a \in \FF_q^\times$ and $c \in \FF_q$ such that
\begin{equation}
\label{eq:degenerate}
    a f(x) = h(x)^p - h(x) + c.
\end{equation}
\end{lemma}
\begin{proof}
Note that if \eqref{eq:degenerate} holds then $h(x)$ is also a Laurent polynomial.
Indeed, write $h(x) = a(x)/b(x)$ with $a(x)$ and $b(x)$ coprime, suppose \eqref{eq:degenerate} holds and note that $h(x)^p - h(x) = (a(x)^p - a(x)b(x)^{p-1})/b(x)^p$ must be a Laurent polynomial.
But if $b(x)^p$ has a non-zero root $z$, then so does $b(x)$, and so does $a(x)$.
Since $a(x)$ and $b(x)$ are coprime, they cannot share a root, and thus $h(x)$ is Laurent.
Write $h(x) = \sum_{s \in S} a_s x^s$ with non-zero $a_s$ and let $n = \max S$ and $m = \min S$.
If $n > 0$, then top term in $h(x)^p - h(x)$ is $a_{s}^{pn} x^{pn}$ which not in $a f(x) - c$ by hypothesis, so $n \leq 0$.
Finally, if $m < 0$, then the lowest term in $h(x)^p - h(x)$ is $a_m x^{pm}$ which is not in $a f(x) - c$ by hypothesis.
\end{proof}

As a corollary of \Cref{thm:cochrane-pinner} and \Cref{lem:degenerate}, we have.

\begin{corollary}
\label{cor:cochrane-pinner}
Let $\FF_q$ be a finite field of characteristic $p$ and let $r_1, \dotsc, r_t$ be distinct non-zero integers, all coprime to $p$.
Write $c \defined \max\set{0,r_1, \dotsc, r_t} + \max\set{0,-r_1, \dotsc, -r_t}$.
Then the sum
\begin{equation*}
    \cS(a_1, \dotsc, a_t, b) \defined \sum_{x \in \FF_q^\times} \psi(a_1 x^{r_1} + \dotsb + a_t x^{r_t}) \chi_b(x),
\end{equation*}
satisfy $\cS(0, \dotsc, 0, 1) = q-1$; $\cS(0, \dotsc, 0, b) = 0$ for all $b \neq 1$; and for each tuple $(a_1, \dotsc, a_t, b) \neq (0, \dotsc, 0, 1)$, we have $\abs[\big]{\cS(a_1, \dotsc, a_t, b)} \leq x q^{1/2}$.
\end{corollary}

With this at hand, we obtain the following general mixing estimate.

\begin{proposition}
\label{prop:general-mixing}
Let $\FF_q$ be a finite field with characteristic $p$ and let $r_1, \dotsc, r_t$ be distinct non-zero integers, all coprime to $p$.
Define $c \defined \max\set{0,r_1, \dotsc, r_t} + \max\set{0,-r_1, \dotsc, -r_t}$.
Then for any sets $A_1, \dotsc, A_t, B \subseteq \FF_q^\times$ with densities $\alpha_i \defined \card{A_i}/q$ and $\beta \defined \card{B}/(q-1)$, define
\begin{equation*}
    X \defined \card[\Big]{\set[\big]{x \in \FF_q \st x^{r_1} \in A_1, \dotsc, x^{r_t} \in A_t , x \in B}},
\end{equation*}
and let $\delta = \alpha_1 \dotsb \alpha_t \beta$ be the product of all the densities.
Then we have
\begin{equation*}
    \abs[\big]{X -  \delta q} \leq c q^{1/2} \norm{\ind_{A_1}}_{\wiener(\FF_q)} \dotsb \norm{\ind_{A_t}}_{\wiener(\FF_q)} \norm{\ind_{B}}_{\wiener(\FF_q^\times)}.
\end{equation*}
\end{proposition}
\begin{proof}
We expand $\ind_{A_1}, \dotsc, \ind_{A_t}$ into a Fourier decomposition in $(\FF_q,+)$ and $\ind_{B}$ into a Fourier decomposition in $(\FF_q^\times,\times)$.
We obtain
\begin{align*}
    X &= \sum_{x \in \FF_q^\times} \ind_{A_1}(x^{r_1}) \dotsb \ind_{A_t}(x^{r_t}) \ind_{B}(x) = \sum_{\substack{a_1, \dotsc, a_t \in \FF_q \\ b \in \FF_q^\times}} \fourier{\ind_{A_1}}(\psi_{a_1}) \dotsb \fourier{\ind_{A_t}}(\psi_{a_t}) \fourier{\ind_B}(\chi_b)\cS(a_1, \dotsc, a_t, b).
\end{align*}
Since $\fourier{\ind_{A_1}}(\psi_{0}) \dotsb \fourier{\ind_{A_t}}(\psi_{0}) \fourier{\ind_B}(\chi_1) = \alpha_1 \dotsb \alpha_t \beta = \delta$, the contribution of the term $(0,\dots,0,1)$ is $\delta \cS(0,\dotsc,0,1) = \delta(q-1)$.
Since $\cS(0, \dotsc, 0, b) = 0$ for $b \neq 1$, we have
\begin{align*}
    \abs[\big]{X - \delta q} 
    &\leq \delta + \!\!\sum_{(a_1, \dotsc, a_t) \neq (0, \dotsc, 0),  b \in \FF_q^\times}\!\! \abs{\fourier{\ind_{A_1}}(\psi_{a_1})} \dotsb \abs{\fourier{\ind_{A_t}}(\psi_{a_t})} \abs{\fourier{\ind_B}(\chi_b)} \abs{\cS(a_1, \dotsc, a_t, b)} \\
    &\leq \delta + c q^{1/2} \!\!\sum_{(a_1, \dotsc, a_t) \neq (0, \dotsc, 0)} \!\! \abs{\fourier{\ind_{A_1}}(\psi_{a_1})} \dotsb \abs{\fourier{\ind_{A_t}}(\psi_{a_t})}  \norm{\ind_{B}}_{\wiener(\FF_q^\times)}\\
    &\leq \delta + c q^{1/2}  \p[\big]{\norm{\ind_{A_1}}_{\wiener(\FF_q)} \dotsb \norm{\ind_{A_t}}_{\wiener(\FF_q)} - \alpha_1 \dotsb \alpha_t} \norm{\ind_{B}}_{\wiener(\FF_q^\times)} \\
    &\leq c q^{1/2} \norm{\ind_{A_1}}_{\wiener(\FF_q)} \dotsb \norm{\ind_{A_t}}_{\wiener(\FF_q)} \norm{\ind_{B}}_{\wiener(\FF_q^\times)}.
\end{align*}
Here we used that $\beta \leq \norm{\ind_{B}}_{\wiener(\FF_q^\times)}$, so $\delta \leq c q^{1/2} \alpha_1 \dotsb \alpha_t \norm{\ind_{B}}_{\wiener(\FF_q^\times)}$.
\end{proof}


\subsection{Supersaturation}

We will make use of a few supersaturation results.
The most general one follows from a version of the Green~\cite{Green2005-nn} removal lemma. 

\begin{lemma}
\label{lem:green-removal}
For every $s \geq 3$ and $\eps > 0$, there is $\delta = \delta_0(s, \eps) > 0$ such that the following holds for all finite commutative groups $G$.
If $A_1, \dotsc, A_s \subseteq G$ are subsets with $(\ind_{A_1} \conv \dotsb \conv \ind_{A_s})(0) \leq \delta$, then for all $1 \leq i \leq s$ there is a subsets $A_i' \subseteq A_i$ with $\card{A_i \setminus A_i'} \leq \eps \card{G}$ such that $(\ind_{A_1'} \conv \dotsb \conv \ind_{A_s'})(0) = 0$.
\end{lemma}

We will use the removal lemma in the following supersaturation guise for finite fields.

\begin{corollary}
\label{cor:green-removal}
For every $L \in \Lin_{\geq 3}$ and $\eps > 0$, there is $\delta = \delta_1(L,\eps) > 0$ such that the following holds.
If $L$ is non-degenerate in $\FF_q$ then every subset $A \subseteq \FF_q$ with $m_L(A) \leq \delta$ satisfy
\begin{equation*}
    \card{A} /q  \leq d_L(\FF_q) + \eps.
\end{equation*}
\end{corollary}
\begin{proof}
Write $L(x_1, \dotsc, x_s) = a_1 x_1 + \dotsb + a_s x_s$ and define the sets $A_1 \defined a_1 \dilation A, \dotsc, A_s \defined a_s \dilation A$.
Since $L$ is non-degenerate in $\FF_q$, the maps $x \mapsto a_i x$ are automorphisms.
Note that $m_L(A) = (\ind_{A_1} \conv \dotsb \conv \ind_{A_s})(0) \leq \delta$.
Choose $\delta_1(L,\eps) = \delta_0(s,\eps/s)$, so by \Cref{lem:green-removal}, we have subsets $Y_i \subseteq A_i$ with $\card{Y_i} \leq \eps q / s$ and with $(\ind_{A_1 \setminus Y_1} \conv \dotsb \conv \ind_{A_s \setminus Y_s})(0) = 0$.
Consider the set
\begin{equation*}
    A' = \set[\big]{x \in \FF_q \st a_1 x \in A_i \setminus Y_1, \dotsc, a_s x \in A_s \setminus Y_s },
\end{equation*}
and note that $A' = A \setminus (a_1^{-1} Y_1 \cup \dotsb \cup a_s^{-1} Y_s)$ satisfy $\card{A'} \geq \card{A} - \card{Y_1} - \dotsb - \card{Y_s} \geq \card{A} - \eps q$ and
\begin{equation*}
    m_L(A') = (\ind_{a_1 \dilation A'} \conv \dotsb \conv \ind_{a_s \dilation A'}) \leq (\ind_{A_1 \setminus Y_1} \conv \dotsb \conv \ind_{A_s \setminus Y_s})(0) = 0,
\end{equation*}
thus $\card{A'}/q \leq d_L(\FF_q)$.
Finally, this gives $\card{A} /  \leq d_L(\FF_q) + \eps$
\end{proof}

The dependency between $\eps$ and $\delta$ given by the removal lemmas above are very poor.
The following result of Balasubramanian and Prakash~\cite[Theorem 2]{Balasubramanian2005-xp} gives explicit quantitative bounds for the supersaturation for sum-free sets.

\begin{proposition}
\label{prop:supersaturation-sumfree}
There is a constant $\delta_0 > 0$ such that for every $0 \leq \delta \leq \delta_0$ and every subset $A \subseteq G$ if a finite commutative group $G$ the following holds.
If $(\ind_A \conv \ind_A \conv \ind_{-A})(0) \leq \delta$, then
\begin{equation*}
    \card{A}/\card{G} \leq d_\Sigma(G) + O(\delta^{1/3}).
\end{equation*}
\end{proposition}

The same result with $O(\delta^{1/5})$ in place of $O(\delta^{1/3})$ is due to Green and Ruzsa~\cite[Proposition 2.2]{Green2005-no}.


\section{Constructions}
\label{sec:constructions}

In this section, we show that there are $L$-free subsets of $\FF_q$ that are close to being maximal density and at the same time, have small Wiener norm.
In the prime field case, we make use of the following result of Candela and Sisask~\cite{Candela2011-cg}.

\begin{theorem}
\label{thm:candela-sisask}
If $L \in \Lin_{\geq 3}$ is not translation invariant, then $\lim_{p \to \infty} d_L(\FF_p) = d_L(\TT) > 0$.
\qed
\end{theorem}

From this, we can deduce that the same limit is attained for $d_L(\FF_p^d)$ when $d$ is fixed.

\begin{lemma}
\label{lem:approximation}
For every $L \in \Lin$ and every $1 \leq h \leq d$, we have
\begin{equation*}
    0 \leq d_L(\FF_p^d) - d_L(\FF_p^h) \leq 1/p^h.
\end{equation*}
\end{lemma}
\begin{proof}
If $Y \subseteq \FF_p^h$ is $L$-free with $\card{Y}/p^h = d_L(\FF_p^h)$ then $A = Y \times \FF_p^{d-h} \subseteq \FF_p^d$ is also $L$-free and has $\card{A}/p^d = d_L(\FF_p^h)$.
Therefore, $d_L(\FF_p^d) \geq d_L(\FF_p^h)$.
For the upper bound, let $A \subseteq \FF_p^d$ be $L$-free with $\card{A} = d_L(\FF_p^d) p^d$ and note that $0 \notin A$.
Write $\Gr_h(\FF_p^d)$ for the set of $h$-dimensional subspaces of $\FF_p^d$.
For each $V \in \Gr_h(\FF_p^d)$, we have that $A \cap V$ is $L$-free and $\card{A \cap V}/p^h \leq d_L(V) = d_L(\FF_p^h)$.
Each non-zero $x \in \FF_p^d$ is in precisely $\card{\Gr_{h-1}(\FF_p^{d-1})}$ such subspace $V$, thus
\begin{equation*}
    d_L(\FF_p^d) = \frac{\card{A}}{p^d} = \!\!\sum_{V \in \Gr_h(\FF_p^d)}\!\! \frac{\card{A \cap V}}{\card{\Gr_{h-1}(\FF_p^{d-1})} p^d} \leq \frac{\card{\Gr_h(\FF_p^d)} d_L(\FF_p^h) p^h}{\card{\Gr_{h-1}(\FF_p^{d-1})} p^d} = \p[\big]{\frac{1 - 1/p^d}{1 - 1/p^h}} d_L(\FF_p^h),
\end{equation*}
where we used that $\card{\Gr_{h}(\FF_p^{d})}/\card{\Gr_{h-1}(\FF_p^{d-1})} = (p^d - 1)/(p^h - 1)$.
Since $d_L(\FF_p^h) \leq 1 - 1/p^h$, we have
\begin{equation*}
    d_L(\FF_p^d) - d_L(\FF_p^h) \leq \p[\Big]{\frac{1 - 1/p^d}{1 - 1/p^h} - 1}d_L(\FF_p^h) = \p[\big]{\frac{1/p^h - 1/p^d}{1 - 1/p^h}} d_L(\FF_p^h) \leq \frac{1}{p^h} - \frac{1}{p^d} \leq \frac{1}{p^h}. \qedhere
\end{equation*}
\end{proof}

As a consequence, Candela and Sisask limit theorem holds more generally.

\begin{corollary}
\label{cor:candela-sisask-d}
If $L \in \Lin_{\geq 3}$, then $d_L(\FF_p^d) = d_L(\TT) + o_p(1)$ for any $d = d(p) \geq 1$ as $p \to \infty$.
\end{corollary}
\begin{proof}
By \Cref{lem:approximation}, we have $d_L(\FF_p^d) = d_L(\FF_p) + o_p(1)$ for any $d = d(p) \geq 1$ and $p \to \infty$, while \Cref{thm:candela-sisask} gives $d_L(\FF_p) = d_L(\TT) + o_p(1)$.
\end{proof}

We now establish that there are near maximisers for $d_L(\FF_p^d)$ with small Wiener norm.

\begin{proposition}
\label{prop:torus-discretisation}
For every $L \in \Lin_{\geq 3}$, and for every $\eps > 0$, there are constants $p_0 = p_0(L, \eps) > 0$ and $C = C(L, \eps) > 0$ such that the following holds.
For every prime $p \geq p_0$ and every $d \geq 1$, there is a $L$-free subset $A \subseteq \FF_p^d$ with $\card{A}/p^d \geq d_L(\FF_p^d) - \eps$ and with $\norm{\ind_A}_\wiener \leq C \log p$.
\end{proposition}
\begin{proof}
If $L$ is translation invariant, then $d_L(\FF_p^d) = o(1)$ as $p \to \infty$, so we can take $A = \emptyset$.
Assume now $L$ is not translation invariant.
The first claim is that there is a $L$-free subset $B \subseteq \TT$ that is a finite union of arcs and with Haar measure $\mu(B) \geq d_L(\TT) - \eps/2$.
We start with a subset $X \subseteq \TT$ that is $L$-free subset and with Haar measure $\mu(X) \geq d_L(\TT) - \eps/4$, given by \Cref{thm:candela-sisask}.
By the inner regularity of the Haar measure, there is a compact subset $K \subseteq X$ with $\mu(K) \geq d_L(\TT) - \eps/2$.
Note that the function $x \mapsto \abs{L(x)}$ is a continuous function $K^s \to \RR$ defined on a compact set $K^s$, and thus, the minimum value $\delta$ of $L(x)$ is attained at a point $x \in K^s$.
Since $K \subseteq X$ is $L$-free, we have $\delta > 0$.
Therefore, there is $\delta' > 0$ such that $B = \cup_{x \in K}B(x,\delta')$, the $\delta'$-neighbourhood of $K$, is also $L$-free.
Indeed, if $L(x_1, \dotsc, x_s) = a_1 x_1 + \dotsc + a_s x_s$, then we can take $\delta' < \delta / 2 (\abs{a_1} + \dotsb + \abs{a_s})$.
Since $K$ is compact, $B$ is a finite union of intervals.
We write $B = I_1 \sqcup \dotsb \sqcup I_m$ for the union of $m = m(L, \eps)$ disjoint intervals and note that $\mu(B) \geq \mu(K) \geq d_L(\TT) - \eps/2$.

Identify $\TT$ with the unit circle $\SS \subseteq \CC$, for every prime $p$, consider the set $Y$ defined as
\begin{equation*}
    Y \defined \set{ y \in \FF_p \st \exp(2 \pi i y / p) \in B}.
\end{equation*}
Since the map $x \mapsto \exp(2 \pi i x / p)$ is a homomorphism, the set $Y$ is also $L$-free.
Now observe that $\lim_{p \to \infty} \card{Y}/p = \mu(B)$ and $d_L(\FF_p^d) = d_L(\TT) + o_p(1)$ by \Cref{cor:candela-sisask-d}, so if $p \geq p_0 = p_0(L, \eps)$ we have $\card{Y} \geq (d_L(\FF_p^d) - \eps)p$. 
Since $Y$ is disjoint union of $m$ arithmetic progressions in $\FF_p$, the triangle inequality and \Cref{prop:wienerinterval} gives
\begin{equation*}
    \norm{\ind_{Y}}_\wiener \leq \sum_{i = 1}^m \norm{\ind_{I_i}}_\wiener \leq m \cI(\FF_p) \leq 2 m \log p.
\end{equation*}
Take $A = Y \times \FF_{p}^{d-1} \subseteq \FF_p^d$ and note that $A$ is $L$-free, satisfy $\card{A}/p^d \geq d_L(\FF_p^d) - \eps$, and since $\ind_A = \ind_Y \otimes \ind_{\FF_p^{d-1}}$,  we have $\norm{\ind_A}_\wiener \leq \norm{\ind_Y}_\wiener \leq 2 m \log p$ from \eqref{eq:wiener-product}.
\end{proof}

The corresponding result in the fixed torsion regime is easier to prove directly from \Cref{lem:approximation}.

\begin{corollary}
\label{cor:small-wiener-coset}
For every $L \in \Lin_{\geq 3}$ and let $1 \leq h \leq d$, there is an $L$-free subset $A \subseteq \FF_p^d$ with $\card{A}/p^d \geq d_L(\FF_p^d) - 1/p^h$ and $\norm{\ind_A}_\wiener \leq p^{h/2}$.
\end{corollary}
\begin{proof}
By \Cref{lem:approximation}, we have $d_L(\FF_p^h)\geq d_L(\FF_p^d) - 1/p^h$. 
Let $B \subseteq \FF_p^h$ be an $L$-free set that attains density $d_L(\FF_p^h)$.
If we take $A = B \times \FF_p^{d-h} \subseteq \FF_q^d$, then we have $\card{A} = \card{B} p^{d-h} = d_L(\FF_p^h) p^d$.
Finally, from \eqref{eq:wiener-product} and \eqref{eq:wiener-cardinality}, we have
\begin{equation*}
    \norm{\ind_A}_{\wiener(\FF_q)} = \norm{\ind_B \otimes \ind_{\FF_p^{d-h}}}_{\wiener(\FF_p^h \times \FF_p^{d-h})} = \norm{\ind_B}_{\wiener(\FF_p^h)} \leq \card{B}^{1/2} \leq p^{h/2}. \qedhere
\end{equation*}
\end{proof}

While we will make no use of the following result, it is conceptually interesting to see that the densities $d_L(\FF_p^d)$ also converge to the density of a limiting object when $p$ is fixed and $d \to \infty$.

\begin{proposition}
\label{prop:limit-fix-char}
Fix a prime $p$ and $L \in \Lin_{\geq 3}$.
Write $\FF_p^\infty = \prod_{i \in \NN} \FF_p$ for the direct product compact commutative group with product topology and product Haar probability measure $\mu$.
Then we have $\lim_{d \to \infty} d_L(\FF_p^d) = d_L(\FF_p^\infty)$ and $0 \leq d_L(\FF_p^\infty) - d_L(\FF_p^h) \leq 1/p^h$ for all $h \geq 1$.
\end{proposition}
\begin{proof}
Since $d_L(\FF_p^d)$ is non-decreasing and bounded, it converges to the limit $\gamma_{p,L} \defined \sup_{d \geq 1} d_L(\FF_p^d)$.
Taking $d \to \infty$ in \Cref{lem:approximation} gives the inequality, so it suffices to show that $d_L(\FF_p^\infty) = \gamma_{p,L}$.

Denote by $\pi_d \from F_p^\infty \to \FF_p^d$ the projection to the first $d$ coordinates.
For each $d \geq 1$ let $B \subseteq \FF_p^d$ be an $L$-free set with $\card{B} = d_L(\FF_p^d) p^d$ and note that the cylinder $A = \set{x \in \FF_p^\infty \st \pi_d(x) \in B}$ is $L$-free, measurable, and has $\mu(A) = d_L(\FF_p^d)$.
This holds for all $d \geq 1$, so $d_L(\FF_p^\infty) \geq \gamma_{p,L}$.

Now let $B \subseteq \FF_p^\infty$ be measurable $L$-free subset with $\mu(B) \geq d_L(\FF_p^\infty) - \eps/2$.
There is $K \subseteq B$ compact with $\mu(K) \geq \mu(B) - \eps/2 \geq d_L(\FF_p^\infty) - \eps$ from the inner regularity Haar probability measure on compact metric spaces.
We now argue that for all sufficiently large $d$, we have that $\pi_d(K) \subseteq \FF_p^d$ is $L$-free.
Indeed, the map $L \from K_p^s \to K_p$ defined as $x \mapsto L(x)$ is continuous and $K^s \subseteq K_p^s$ is a compact subset.
Therefore, $L(K^s)$ is a compact and since $K$ is $L$-free, we have $0 \notin L(K^s)$.
Thus, $K_p \setminus L(K^s)$ is an open neighbourhood of $0$, and hence, it must contain $H_d \defined \ker \pi_d$ for all sufficiently large $d$.
For such $d$, we have then that $A_d \defined \pi_d(K) \subseteq \FF_p^d$ is $L$-free, and that $\pi_d^{-1}(A_d) \subseteq \FF_p^\infty$ is $L$-free and contain $K$.
Therefore, for all $\eps > 0$, we get
\begin{equation*}
    \gamma_{p,L} \geq d_L(\FF_p^d) \geq \card{A_d}/p^d = \mu(\pi_d^{-1}(A_d)) \geq \mu(K) \geq d_L(\FF_p^\infty) - \eps. \qedhere
\end{equation*}
\end{proof}

In any case, the main technical result we need for our constructions is the existence of almost maximisers with small Wiener norm for all fields.

\begin{corollary}
\label{cor:small-wiener-near-maximiser}
For every $L \in \Lin_{\geq 3}$ and $\eps > 0$, there is $q_0 = q_0(L, \eps) > 0$ such that the following holds for $q \geq q_0$: there is an $L$-free set $A \subseteq \FF_q$ with $\card{A}/q \geq d_L(\FF_q) - 
\eps$ and $\norm{\ind_A}_\wiener \leq (\log q)^2$.
\end{corollary}
\begin{proof}
Fix $\eps > 0$.
From \Cref{prop:torus-discretisation}, there is an $L$-free subset $A \subseteq \FF_q$ with $\card{A}/q \geq d_L(\FF_q) - \eps$ and $\norm{\ind_A}_\wiener \leq C \log p$ if $q = p^d$ and $p \geq p_0 = p_0(L,\eps) > 0$ and where $C = C(L,\eps) > 0$.
Thus if $p \geq p_0$, we have $\norm{\ind_A}_\wiener \leq C(L,\eps) \log p$, which is less than $(\log q)^2$ for large enough $q$.

Now assume $p \leq p_0$.
From \Cref{cor:small-wiener-coset}, there is an $L$-free subset $B \subseteq \FF_q$ with $\card{B}/q \geq d_L(\FF_q) - 1/p^h$ and $\norm{\ind_B}_\wiener \leq p^{h/2}$ for each $1 \leq h \leq d$.
Take $h = \ceil{\log_p(1/\eps)}$ and note $d \geq h$ for $q \geq p_0^h$.
We have then $\card{B}/q \geq d_L(\FF_q) - \eps$ and $\norm{\ind_B}_\wiener \leq (p_0/\eps)^{1/2} \leq (\log q)^2$ for large enough $q$.
\end{proof}

Finally, we can prove the lower bound in \Cref{thm:general-independence-fq}, implies also the lower bound on \Cref{thm:density-equation-inverse-fq} and \Cref{thm:density-independence-fq}.

\begin{proposition}
\label{prop:general-construction}
For every $L_1, \dotsc, L_t \in \Lin_{\geq 3}$, distinct integers $r_1, \dotsc, r_t$, and $\eps > 0$, there $q_0 > 0$ such that the following is true.
For every $q \geq q_0$ such that $r_1, \dotsc, r_t$ are non-zero in $\FF_q$, there is a set $A \subseteq \FF_q$ such that $A^{r_i}$ is $L_i$-free for every $1 \leq i \leq t$ with
\begin{equation*}
    \card{A}/q \geq d_{L_1}(\FF_q) \dotsb d_{L_t}(\FF_q)  - \eps.
\end{equation*}
Further, if $t = 2k$, with $r_{2i} = - r_{2i-1}$ and $L_{2i} = L_{2i-1}$ for $1 \leq i \leq k$, then we may assume $A = A^{-1}$.
\end{proposition}
\begin{proof}
From \Cref{cor:small-wiener-near-maximiser}, for every $\eps' > 0$, there are sets $X_1, \dotsc, X_t \subseteq \FF_q$ such that $X_i$ is $L_i$-free, $\card{X_i}/q \geq d_{L_i}(\FF_q) - \eps'$, and $\norm{\ind_{X_i}}_\wiener \leq (\log q)^2$ for all $1 \leq i \leq t$.
Consider the set
\begin{equation*}
    A \defined \set[\big]{x \in \FF_q^\times \st x^{r_1} \in X_1, \dotsc , x^{r_t} \in X_t},
\end{equation*}
and note that for each $1 \leq i \leq r$, we have that $A^{r_i} \subseteq X_{i}$, and hence $A^{r_i}$ is $L_i$-free.
If the additional symmetry condition is given, take $X_{2i} = X_{2i-1}$ for $1 \leq i \leq k$ and note that $A = A^{-1}$.

In any case, we use the mixing estimate \Cref{prop:general-mixing} to obtain
\begin{align*}
    \card{A}/q &\geq (d_{L_1}(\FF_q) - \eps') \dotsb (d_{L_t}(\FF_q) - \eps')  - c q^{-1/2} (\log q)^{2t} \\
    &\geq d_{L_1}(\FF_q) \dotsb d_{L_t}(\FF_q) - (t+1) \eps'
\end{align*}
for large enough $q$.
Choose $\eps' = \eps/(t+1)$ and the result follows.
\end{proof}

We now address the enumeration component of \Cref{thm:general-independence-fq} and consequently of \Cref{thm:counting-independence-fq}. 
The lower bound on the enumeration part of \eqref{eq:general-independence-fq} follows directly.
Pick a maximum density set $A$ in $\cF = \cF_{L_1, \dotsc, L_t}^{(r_1, \dotsc, r_t)}(\FF_q)$, thus all $2^{\card{A}}$ of its subsets are also in $\cF$.
For \Cref{thm:counting-equation-inverse-fq}, note that if $A \subseteq \FF_q^\times$ is self-inverse, then there are at least $2^{(\card{A}-k)/2}$ self-inverse subsets of $A$, where $k = \card{\set{x \in A \st x = x^{-1}}} \leq 2$.

To construct product-free sets for \Cref{thm:general-independence-fq-product}, we need large sum-free subsets of $\ZZ/N\ZZ$ with small Wiener norm.
For sum-free sets, we know more.

\begin{lemma}
\label{lemma:sumfree-construction}
For every finite commutative group $G$, there are maximum-size sum-free subsets $A \subseteq G$ with $\norm{\ind_A}_\wiener \leq 2 \log \card{G}$.
\end{lemma}
\begin{proof}
Note that the constructions of maximum sum-free sets Green and Ruzsa~\cite{Green2005-no} are all of the form $I \times H$ where $I$ is an interval and $H$ is a subgroup of $G$.
The result follows immediately from \eqref{eq:wiener-product} and \Cref{cor:phi-G}.
\end{proof}

We now furnish a construction for \Cref{thm:general-independence-fq-product}.

\begin{proposition}
\label{prop:general-construction-product}
For every $L_1, \dotsc, L_t \in \Lin_{\geq 3}$, distinct integers $r_1, \dotsc, r_t$, and $\eps > 0$, there $q_0 > 0$ such that the following is true.
For every $q \geq q_0$ such that $r_1, \dotsc, r_t$ are non-zero in $\FF_q$, there is a product-free set $A \subseteq \FF_q$ such that $A^{r_i}$ is $L_i$-free for every $1 \leq i \leq t$ with
\begin{equation*}
    \card{A}/q \geq d_{L_1}(\FF_q) \dotsb d_{L_t}(\FF_q)d_\Pi(\FF_q^\times) - \eps.
\end{equation*}
\end{proposition}
\begin{proof}
Choose $\eps' = \eps/(t+1)$ and take $X_1, \dotsc, X_t \subseteq \FF_q$ as in the proof of \Cref{prop:general-construction}.
Since $(\FF_q^\times, \times)$ is isomorphic to $(\ZZ/(q-1)\ZZ)$, \Cref{lemma:sumfree-construction} gives us a product-free subset $B \subseteq \FF_q^\times$ with $\card{B}/(q-1) = d_\Pi(\FF_q^\times)$ and $\norm{\ind_{B}}_{\wiener(\FF_q^\times)} \leq 2\log q$.
Then the set
\begin{equation*}
    A \defined \set[\big]{x \in \FF_q^\times \st x^{r_1} \in X_1, \dotsc , x^{r_t} \in X_t, x \in B}
\end{equation*}
is product-free and $A^{r_i}$ is $L_i$-free for $1 \leq i \leq t$.
By \Cref{prop:general-mixing}, for large enough $q$ we have
\begin{align*}
    \card{A}/q &\geq (d_{L_1}(\FF_q) - \eps') \dotsb (d_{L_t}(\FF_q) - \eps')d_\Pi(\FF_q^\times)  - c q^{-1/2} (\log q)^{2t + 1} \\
    &\geq d_{L_1}(\FF_q) \dotsb d_{L_t}(\FF_q)d_\Pi(\FF_q^\times) - \eps. \qedhere
\end{align*}
\end{proof}


\section{Wrapper Theorem}
\label{sec:wrapper}

Our discussion of wrappers follows closely to the one by Semchankau~\cite{Semchankau2025-nt}.

We fix the following standard partition of the unit circle $\SS \subseteq \CC$.
For each positive integer $K$, we consider the partition $\SS = \SS_{0,K} \sqcup \dotsb \sqcup \SS_{K-1,K}$, where we have
\begin{equation*}
    \SS_{j,K} \defined \set[\big]{\exp(2\pi \theta / K) \st j - 1/2 \leq \theta < j + 1/2 } \subseteq \SS.
\end{equation*}
We the maps $m_K \from \SS \to \SS$ that sends each arc to it's midpoint.
That is, $m_K$ is the map defined by $m_K(\SS_{j,K}) = \exp(2\pi j /K)$.
Note that if $K \geq 2$ and $z \in \SS$, then $\abs{z - m_K(z)} \leq \pi/K$.

\begin{definition}
Let $G$ be a finite commutative group and $K \geq 1$ be an integer.
Let $\Lambda \subseteq \fourier{G}$ be a set of characters and let $\sigma \from \Lambda \to \ZZ/K\ZZ$ be any function.
Then the block $B(\Lambda,\sigma) \subseteq G$ is the set
\begin{equation*}
   B(\Lambda, \sigma) \defined \set[\big]{x \in G \st \chi(x) \in \SS_{\sigma(\chi),K}, \text{for all $\chi \in \Lambda$} }. 
\end{equation*}
A block of this form is said to be a $(\Lambda,K)$-block.
A \emph{$(\Lambda, K)$-wrapper} is a disjoint union of $(\Lambda, K)$-blocks, that is, it is a set of the form
\begin{equation*}
    W = \bigcup \set[\big]{B(\Lambda, \sigma) \st \sigma \in \cC}
\end{equation*}
for some collection $\cC$ of maps $\Lambda \to \ZZ/K\ZZ$.
\end{definition}

If $\sigma \from \Lambda \to \ZZ/K\ZZ$ is the constant zero map, then the block $B(\Lambda, \sigma)$ is a Bohr set.
A block more generally an \emph{inhomogeneous Bohr set} where for each prescribed frequency $\chi$ in $\Lambda$ is centred at some point $m_K(\SS_{\sigma(\chi),K})$.

We remark that Semchankau~\cite{Semchankau2025-nt} defines $(\eps,d)$-blocks and $(\eps,d)$-wrappers rather than our $(\Lambda, K)$-blocks, where $d = \card{\Lambda}$ is the \emph{dimension} and where $\eps \defined 1/K$ is the \emph{granularity} and must be taken so that $1/\eps$ is an integer.
We prefer making the frequency set explicit.

One of the main advantages of working with wrappers is that their Wiener norm is controlled.

\begin{proposition}
\label{prop:wiener-wrapper}
For every $(\Lambda,K)$-wrapper $W \subseteq G$, we have $\norm{\ind_W}_\wiener \leq \p{K{\cB(G)}}^{\card{\Lambda}}$.
\end{proposition}
\begin{proof}
We first show that for a block $B = B(\Lambda, \sigma)$, we have
\begin{align*}
    \norm{\ind_B}_\wiener = \norm[\Big]{\prod_{\chi \in \Lambda}\ind_{B(\set{\chi}, \sigma)}}_\wiener \leq \prod_{\chi \in \Lambda} \norm[\big]{\ind_{\chi^{-1}(\SS_{\sigma(\chi),K})}}_\wiener \leq \cB(G)^{\card{\Lambda}},
\end{align*}
from \eqref{eq:algebra-norm}.
Therefore, if $W = \bigcup \set{B(\Lambda, \sigma) \st \sigma \in \cC}$ is a $(\Lambda, K)$-wrapper.
Note that $\card{\cC} \leq K^{\card{\Lambda}}$, so the triangle inequality gives
\begin{equation*}
    \norm{\ind_W}_\wiener \leq \sum_{\sigma \in \cC} \norm{\ind_{B(\Lambda,\sigma)}}_\wiener \leq \card{\cC} \cB(G)^{\card{\Lambda}} \leq (K \cB(G))^{\card{\Lambda}}. \qedhere
\end{equation*}
\end{proof}

For our enumeration results, the following estimate on the number of wrappers is used.

\begin{lemma}
\label{lemma:wrapper-count}
The number of $(\Lambda,K)$-wrappers in $G$ with $\card{\Lambda} \leq d$ is at most $2^{K^d + O(d \log \card{G})}$.
\end{lemma}
\begin{proof}
To fully describe a wrapper, we need to choose a set of characters $\Lambda \subseteq \fourier{G}$ with $\card{\Lambda} \leq d$ and we need to choose a subset of maps $\Lambda \to \ZZ / K \ZZ$.
Thus, the number of wrappers is bound by
\begin{equation*}
    \sum_{k = 0}^{d}\binom{\card{\fourier{G}}}{k}2^{K^k} \leq 2^{K^d} \p[\Big]{\frac{e \card{G}}{d}}^d \leq 2^{K^d + O(d \log \card{G} )}. \qedhere
\end{equation*}
\end{proof}

The main result technical result of this section is the following theorem, which is a more ergonomic modification of Semchankau~\cite[Wrapping Lemma]{Semchankau2025-nt} and the version of Chow, Lim, and Mudgal~\cite[Theorem 3.1]{Chow2025-hs}.

\begin{theorem}
\label{thm:polywrapper}
Let $s \geq 3$ be an integer and $0 < \tau, \rho \leq 1$, $0 \leq \eps \leq 1$ be real parameters satisfying
\begin{equation}
\label{eq:wrapper-parameters}
    \tau \rho^s \geq 2 (4s)^s \eps.
\end{equation}
Then for any subsets $A_1, \dotsc, A_s \subseteq G$ of a finite commutative group $G$ satisfying
\begin{equation}
\label{eq:wrapper-convolution-condition}
    (\ind_{A_1} \conv \dotsb \conv \ind_{A_s} \conv \kappa)(a) \leq \eps,
\end{equation}
where $a \in G$ and $\kappa \from G \to \RR_{\geq 0}$ is a kernel function satisfying $\norm{\kappa}_1 = 1$, the following must hold.
There exists an integer
\begin{equation}
\label{eq:wrapper-bound-k}
    K \ll_s \tau^{-1} \rho^{1 - s},
\end{equation}
and for every $1 \leq i \leq s$, exceptional sets $Y_i \subseteq A_i$ with $\card{Y_i} \leq \rho \card{G}$, frequency sets $\Lambda_i \subseteq \fourier{G}$ with 
\begin{equation}
\label{eq:wrapper-bound-lambda}
    \card{\Lambda_i} \ll_s \tau^{-2}\rho^{2 - 2s}\log(1/\rho \tau),
\end{equation}
and $(\Lambda_i, K)$-wrappers $W_i \subseteq G$ with $A_i \setminus Y_i \subseteq W_i$, and with the property that
\begin{equation}
\label{eq:wrapper-convolution-conclusion}
    (\ind_{W_1} \conv \dotsb \conv \ind_{W_s} \conv \kappa)(a) \leq \tau.
\end{equation}

Moreover, one can additionally prescribe a partition $\Sigma$ of $\set{1, \dotsc, s}$, a set $T_\sigma \subseteq G$ for each $\sigma \in \Sigma$, and automorphisms $\varphi_i \in \Aut(G)$ such that $A_i = \varphi_i(T_\sigma)$ for all $i \in \sigma$, $\sigma \in \Sigma$.
Then, we may assume that $W_i = (\varphi_i \compo \varphi_j^{-1})(W_j)$ and $Y_i = (\varphi_i \compo \varphi_j^{-1})(Y_j)$ for all $i,j \in \sigma$, for all $\sigma \in \Sigma$.
\end{theorem}

We have a few remarks about the statement of \Cref{thm:polywrapper} before we go to the proof:

\begin{enumerate}
    \item Note that we can ignore $\kappa$ in the convolution by setting $\kappa = \card{G}\ind_{\set{0}}$, so condition \eqref{eq:wrapper-convolution-condition} is equivalent to $(\ind_{A_1} \conv \dotsb \conv \ind_{A_s})(a) \leq \eps$, and \eqref{eq:wrapper-convolution-conclusion} is equivalent to $(\ind_{W_1} \conv \dotsb \conv \ind_{W_s})(a) \leq \tau$.
    \item If we want to use a kernel $\kappa \from G \to \RR_{\geq 0}$ with a general value of $\mu = \norm{\kappa}_1 > 0$, then we can use the kernel $\kappa/\mu$ and change the constant $\eps$ to $\eps/\mu$ in \eqref{eq:wrapper-convolution-condition}.
    \item If we take $\kappa = \ind_{A_{s+1}} \conv \dotsb \conv \ind_{A_{k}}$, then this result can be seen as a \emph{partial wrapper theorem}, where we find wrappers for some but not all sets in the convolution. 
    \item It would straightforward to write non-uniform and explicit bounds for $K$ and $\card{\Lambda_i}$ that depend of $i$ at the expense of simplicity and ease of use.
    Moreover, the dependencies on $s$ of the constants implicit in $\ll_s$ in \eqref{eq:wrapper-bound-k} and \eqref{eq:wrapper-bound-lambda} are not optimal in the proof.
    \item We may always assume that we have a partition $\Sigma$ in \Cref{thm:polywrapper}, since if none is given, we can take a partition of $\set{1, \dotsc, s}$ into singletons and each $\varphi_i$ being the identity map.
    \item \Cref{thm:polywrapper} is also applicable when $\eps = 0$.
    Condition \eqref{eq:wrapper-parameters} becomes vacuous.
    \item Comparing with Chow, Lim, and Mudgal~\cite[Theorem 3.1]{Chow2025-hs} and with Semchankau~\cite[Wrapping Lemma]{Semchankau2025-nt}, our version requires no additional shift from in the convolution: namely, both convolutions \eqref{eq:wrapper-convolution-condition} and \eqref{eq:wrapper-convolution-conclusion} are evaluated at the same point $a \in A$.
    Technically, both their versions assume that $\card{A_i} \gg \card{G}$, and while we do not explicitly assume this, the conclusion \eqref{eq:wrapper-convolution-conclusion} may not useful without this assumption.
    The version in \cite{Chow2025-hs} also allows for arbitrary rotations of the arc partition $\SS = \SS_{0,K} \sqcup \dotsb \sqcup \SS_{K-1,K}$, which we avoid.
    We also get polynomial dependencies on $\eps$, while their parameters depend exponentially on $\eps$.
    Finally, we encode certain symmetries in the wrappers that will be convenient for applications.
\end{enumerate}

We now present a full proof of \Cref{thm:polywrapper}.
We require the following Fourier sampling result of Croot, Łaba, and Sisask~\cite{Croot2013-vy}.

\begin{proposition}
\label{prop:almost-periodicity}
Let $G$ be a finite commutative group.
Let $0 < \tau < 1$ and $r > 2$ be parameters.
For any non-zero function $f \from G \to \CC$, there are characters $\chi_1, \dotsc, \chi_d \in \fourier{G}$ with $d \leq 20 r / \tau^2$ and coefficients $c_1, \dotsc, c_d \in \CC$ with $\abs{c_j} \leq 1$ such that
\begin{equation*}
    \norm[\Big]{\frac{f}{\norm{f}_{\wiener}} - \frac{1}{d}\sum_{j \in [d]} c_j \chi_j }_r \leq \tau.
\end{equation*}
\end{proposition}

A version that gives the explicit constant of $20$ above is in Ramaré~\cite{Ramare2024-iq}.
Equipped with this, we can prove the following decomposition theorem due to Semchankau~\cite[Wrapper Theorem]{Semchankau2025-nt} with explicit constants.

\begin{lemma}
\label{lem:wrapper-decomposition}
Let $f \from G \to \CC$ be a non-zero function with Wiener norm $\norm{f}_\wiener = w$.
Let $0 < \xi < 1/4$ and $0 < \delta < w$ be real numbers and $K \geq 9 w/\delta$ be an integer.
Then there exists a decomposition
\begin{equation*}
    f = f_{\str} + f_{0} + f_{\infty}
\end{equation*}
of $f$ into functions $f_{\str}, f_{0}, f_{\infty} \from G \to \CC$ such that $\norm{f_{\infty}}_\infty \leq \delta$, $\card{\supp{f_{0}}} \leq \xi \card{G}$ and there is a set of characters $\Lambda \subseteq \fourier{G}$ with $\card{\Lambda} \leq (2^{11} w^2 /\delta^2 ) \log_2(1/\xi)$ such that $f_{\str}$ is constant on each $(\Lambda, K)$-block.
\end{lemma}
\begin{proof}
Apply \Cref{prop:almost-periodicity} with $r \defined \log_2(1/\xi)$ and $\tau \defined \delta/9w$.
We obtain a set of characters $\Lambda = \set{\chi_1, \dotsc, \chi_d} \subseteq \fourier{G}$ with $\card{\Lambda} \leq d \leq 20 r / \tau^2 \leq (2^{11} w^2 /\delta^2 ) \log_2(1/\xi)$ and coefficients $c_1, \dotsc, c_d \in \CC$ with $\abs{c_j} \leq 1$ such that if we write
\begin{equation*}
    f_{\Lambda} \defined \frac{w}{d}\sum_{j \in [d]} c_j \chi_j,
\end{equation*}
then $f_{r} \defined f - f_{\Lambda}$ satisfies $\norm{f_{r}}_r \leq w \tau$.
Partition of $G$ into $(\Lambda,K)$-blocks $B(\Lambda, \sigma)$ with $\sigma \from \Lambda \to \ZZ / K\ZZ$.
For each such $\sigma$, choose an element $x_\sigma \in B(\Lambda, \sigma)$.
We define our function $f_{str} \from G \to \CC$ to be constant on each block $B(\Lambda, \sigma)$ with $f_{str}(x) \defined f_{\Lambda}(x_\sigma)$ for $x \in B(\Lambda, \sigma)$.
Define $h_1 = f_{\Lambda} - f_{str}$ and note that for every $x \in B(\Lambda,\sigma)$, we have
\begin{align*}
    \abs{h_1(x)} &= \abs{f_{\Lambda}(x) - f_{\Lambda}(x_\sigma)} = \abs[\bigg]{\frac{w}{d} \sum_{j \in [d]} c_j \p[\big]{\chi_j(x) - \chi_j(x_\sigma)} } \\
    &\leq \frac{w}{d} \sum_{j \in [d]} \abs{c_j}\abs{\chi_j(x) - \chi_j(x_\sigma)} \leq \frac{w}{d}\sum_{j \in [d]} \frac{2\pi}{K} \leq \frac{2\pi w}{K} \leq 2 \pi w \tau,
\end{align*}
where here we used that if, $z, z'$ lie in an arc $A \subseteq \SS \subseteq \CC$ of $\alpha$ radians, then, then $\abs{z - z'} \leq 2 \sin(\alpha/2) \leq \alpha$.
So far we have the decomposition $f = f_{\Lambda} + f_{r} = (f_{str} + h_1) + f_{r}$, and we will now further decompose $f_r$.
Let $S \defined \set{x \in G \st f_r(x) \geq 2 w\tau}$ and define $f_{0} \defined f_{r} \ind_{S}$ and $h_2 \defined f_{r} - f_{0}$.
Estimate the size of $S$ as follows,
\begin{equation*}
    \card{S} (2 w\tau)^r \leq \sum_{x \in S}\abs{f_{r}(x)}^r \leq \card{G}\, \norm{f_{r}}_r^r \leq \card{G} (w \tau)^r,
\end{equation*}
therefore, $\card{\supp f_{0}} \leq \card{S} \leq 2^{-r} \card{G} = \xi \card{G}$.
Also, note that $\norm{h_2(x)}_{\infty} \leq 2 w \tau$.
Set $f_{\infty} \defined h_1 + h_2$, so we have $f = f_{str} + f_0 + f_\infty$ with $\norm{f_{\infty}}_\infty \leq \norm{h_1}_\infty + \norm{h_2}_\infty \leq (2\pi + 2) w \tau \leq 9 w \tau = \delta$.
\end{proof}

Note that $f \from G \to \RR$ is a real function, then we can take $f_{\str}, f_0, f_\infty \from G \to \RR$ real as well.
Using \Cref{lem:wrapper-decomposition}, we can approximate the level set of functions by wrappers.

\begin{corollary}
\label{cor:wrapper-level}
Let $f \from G \to \RR$ be a non-zero function with Wiener norm $\norm{f}_{\wiener} = w$.
Let $\ell_1 \leq \ell_2$, $0 < \xi < 1/4$, and $\delta > 0$ be real numbers, and let $K \geq 9 \max\set{4 w / \delta, 1}$ be an integer.
Then there exists a $(\Lambda, K)$-wrapper $V \subseteq G$ with
$\card{\Lambda} \leq 2^{11} \max\set{16 w^2/\delta^2, 1} \log_2(1/\xi)$, and an exceptional set $Z \subseteq G$ with $\card{Z} \leq \xi \card{G}$, such that:
\begin{equation}
\label{eq:level-set-inner}
    \set[\big]{x \in G \st \ell_1 \leq f(x) \leq \ell_2} \setminus Z \subseteq V,
\end{equation}
\begin{equation}
\label{eq:level-set-outer}
    V \setminus Z \subseteq \set[\big]{x \in G \st \ell_1 - \delta \leq f(x) \leq \ell_2 + \delta}.
\end{equation}
\end{corollary}
\begin{proof}
Apply \Cref{lem:wrapper-decomposition} with the value of $\delta$ replaced by $\min\set{w,\delta/4}$.
We obtain a corresponding decomposition $f = f_{\str} + f_0 + f_\infty$.
Let $Z = \supp f_0$ so $\card{Z} \leq \chi \card{G}$ and set
\begin{equation*}
    V \defined \set[\big]{x \in G \st \ell_1 - \delta/2 \leq f_{\str}(x) \leq \ell_2 + \delta/2}.
\end{equation*}
Since $f_{\str}$ is constant on $(\Lambda,K)$-blocks, if one element $x \in B(\Lambda, \sigma)$ is in $V$, then $B(\Lambda, \sigma) \subseteq V$.
Therefore, $V$ is a $(\Lambda, K)$-wrapper.
Finally, observe that $\norm{f_\infty}_\infty \leq \delta/4$ and that for $x \notin Z$, we have $f(x) = f_{\str}(x) + f_\infty(x)$.
Thus, if $x \notin Z$ and $\ell_1 \leq f(x) \leq \ell_2$, then $\ell_1 - \delta/4 \leq f_{\str}(x) \leq \ell_2 + \delta/4$, and so $x \in V$ and \eqref{eq:level-set-inner} holds.
Similarly, $x \notin Z$ and $f(x) \leq \ell_1 - \delta$ (or $f(x) \geq \ell_2 + \delta$) then $f_{\str}(x) \leq \ell_1 - 3\delta/4$ (or $f_{\str}(x) \geq \ell_2 + 3\delta/4$), so $x \notin V$ and \eqref{eq:level-set-outer} holds.
\end{proof}

We finally can proceed to the proof of the main technical result.

\begin{proof}[Proof of \Cref{thm:polywrapper}]
Define the parameters $t_0 \leq t_1 \leq \dotsb \leq t_s$ and $\xi$ as
\begin{equation*}
    t_i \defined \frac{\tau}{4}\p[\Big]{\frac{\rho}{4s}}^{s-i}, \qquad \xi \defined \min\set[\Big]{\frac{\rho}{2s},\frac{\tau}{2s}},
\end{equation*}
thus we have $t_0 \geq \eps/2$ by hypothesis, $t_s = \tau/4$, and $t_{i-1}/t_i = \rho/4s$.
Also note that $\xi < 1/4$ and that $\log_2(1/\xi) \ll_s \log(1/\rho \tau)$.
Let $K \defined \ceil{36/t_1}$ and note that $K \geq 9 \min\set{1,4/t_i}$ for all $1 \leq i \leq s$, and also that $K \ll_s \tau^{-1}\rho^{1-s}$ as claimed in \eqref{eq:wrapper-bound-k}.

Our first step is to construct a preliminary family of $(\Lambda_i, K)$-wrappers $V_i$ and later we will modify them so they satisfy the prescribed automorphisms.
Indeed, we sequentially build $V_i$, exceptional sets $Z_i$, and auxiliary sets $U_i = V_i \setminus Z_i$ by applying \Cref{cor:wrapper-level} in the order $i = 1, \dotsc, s$.

Assume inductively that $U_1, \dotsc, U_{i-1}$ have been determined.
Define the function $F_i \from G \to \RR_{\geq 0}$ by
\begin{equation*}
    F_i(x) \defined (\ind_{U_1} \conv \dotsb \conv \ind_{U_{i-1}} \conv \ind_{A_{i+1}} \conv \dotsb \conv \ind_{A_s} \conv \kappa)(a - x),
\end{equation*}
and note that $w_i \defined \norm{F_i}_\wiener \leq 1$ by \eqref{eq:wiener-convolution}.
If $F_i = 0$, take $V_i = U_i = G$, $Z_i = \emptyset$, $\Lambda_i = \emptyset$.
Otherwise apply \Cref{cor:wrapper-level} with $\ell_1 = 0$, $\ell_2 = t_i$, $\delta = t_i$, and the constant $\xi$.
We obtain an exceptional set $Z_i \subseteq G$, a $(\Lambda_i, K)$-wrapper $V_i \subseteq G$ and set $U_i \defined V_i \setminus Z_i$.

At the end of this process, we have constructed sets $Z_i$ with $\card{Z_i} \leq \xi \card{G}$, as well as $(\Lambda_i, K)$-wrappers $V_i$ with frequency sets $\Lambda_i$ satisfying \eqref{eq:wrapper-bound-lambda} since
\begin{equation*}
    \card{\Lambda_i} \leq 2^{11} \max\set{16 w_i^2/t_i^2,1} \log_2(1/\xi) \ll_s \tau^{-2}\rho^{2-2s} \log(1/\rho \tau),
\end{equation*}
further, the sets $U_i \defined V_i \setminus Z_i$ are such that $F_i(x) \leq 2 t_i$ for $x \in U_i$, and 
\begin{align*}
    \set{x \in G \st F_i(x) \leq t_i} \setminus Z_i \subseteq V_i.
\end{align*}
For $0 \leq i \leq s$, define
\begin{equation*}
    L_i = (\ind_{U_1} \conv \dotsb \conv \ind_{U_{i-1}} \conv \ind_{U_{i}}  \conv \ind_{A_{i+1}} \conv \dotsb \conv \ind_{A_s} \conv \kappa)(a).
\end{equation*}
Note that $L_i = \inner{\ind_{U_i}}{F_i}$ of $1 \leq i \leq s$, hence
\begin{equation}
\label{eq:Li-bound}
    L_i = \inner{\ind_{U_i}}{F_i} = \frac{1}{\card{G}} \sum_{x \in G} \ind_{U_i}(x) F_i(x) \leq 2 t_i.
\end{equation}
Note that \eqref{eq:Li-bound} also holds for $i = 0$, since $L_0 \leq \eps$ by hypothesis \eqref{eq:wrapper-convolution-condition}, and $\eps \leq 2t_0$ by hypothesis~\eqref{eq:wrapper-parameters}.
In any case, the following bound on the convolution of the $U_i$'s holds
\begin{equation}
\label{eq:Ls-bound}
    L_s = (\ind_{U_1} \conv \dotsb \conv \ind_{U_{s}} \conv \kappa)(a) \leq 2 t_s.
\end{equation}
To turn this into an upper bound on the convolution of the $V_i$'s, we simply replace the indicator $\ind_{U_i}$ by $\ind_{V_i} = \ind_{U_i} + \ind_{Z_i \cap V_i}$ one index $i$ at a time.
Using \eqref{eq:Ls-bound}, together with fact that
\begin{equation*}
    (\ind_{X_1} \conv \dotsb \conv \ind_{X_s} \conv \kappa)(a) \leq (\ind_{X_i} \conv \ind_G \conv \kappa)(a) = \card{X_i}/\card{G}
\end{equation*}
for all $1 \leq i \leq s$ and for sets $X_i$, linearity of the convolution gives
\begin{align*}
    (\ind_{V_1} &\conv \dotsb \conv \ind_{V_s} \conv \kappa)(a) \leq (\ind_{V_1} \conv \dotsb \conv \ind_{V_{s-1}} \conv \ind_{U_s} \conv \kappa)(a) + \card{Z_s}/\card{G} \\
    &\leq \dotsb \leq (\ind_{U_1} \conv \dotsb \conv \ind_{U_s} \conv \kappa)(a) + \p{\card{Z_1} + \dotsb + \card{Z_s}}/\card{G}
    \leq 2 t_s + s \xi \leq \frac{\tau}{2} + \frac{\tau}{2} = \tau.
\end{align*}

Next we show that $\card{A_i \setminus V_i}$ cannot be too large.
Markov's inequality gives the bound
\begin{equation*}
    L_{i-1} \card{G} = \sum_{x \in G} \ind_{A_i}(x) F_i(x) \geq \sum_{x \in G} \ind_{A_i}(x) \ind_{\set{F_i(x) > t_i}} F_i(x) > t_i \card[\big]{\set{x \in A_i \st F_i(x) > t_i}}. 
\end{equation*}
Hence, we have
\begin{align*}
    \card{A_i \setminus V_i} &= \card[\big]{\set{x \in A_i \setminus V_i \st F_i(x) \leq t_i}} + \card[\big]{\set{x \in A_i \setminus V_i \st F_i(x) > t_i}} \\
    &\leq \card{Z_i} + \card[\big]{\set{x \in A_i \st F_i(x) > t_i}} \\
    &\leq \xi \card{G} + \frac{L_{i-1}}{t_i} \card{G} 
    \leq \p[\big]{\xi + \frac{2t_{i-1}}{t_i}}\card{G} \leq \p[\big]{\frac{\rho}{2s} + \frac{\rho}{2s}}\card{G} = \frac{\rho}{s} \card{G}.
\end{align*}

Now we impose the prescribed automorphism relations.
For each $\sigma \in \Sigma$, define
\begin{equation*}
    X_\sigma \defined \bigcap_{j \in \sigma} \varphi_j^{-1}(V_j),
\end{equation*}
and for each $i \in \sigma$, we set $W_i \defined \varphi_i(X_\sigma)$.
Hence, for $i,j \in \sigma$, we have $(\varphi_i \compo \varphi_j^{-1})(W_j) = \varphi_i(X_\sigma) = W_i$ as claimed.
Further, it holds that $W_i \subseteq V_i$ and hence, by monotonicity,
\begin{equation*}
    (\ind_{W_1} \conv \dotsb \conv \ind_{W_s} \conv \kappa)(a) \leq \tau,
\end{equation*}
so \eqref{eq:wrapper-convolution-conclusion} is established.
Define the exceptional set $Y_i \defined A_i \setminus W_i \subseteq A_i$ and note that, if $i \in \sigma$, then
\begin{align*}
    \card{Y_i} &= \card{A_i \setminus W_i} = \card{\varphi_i^{-1}(A_i) \setminus \varphi_i^{-1}(W_i)} \\
    &= \card{T_\sigma \setminus X_\sigma} \leq \sum_{j \in \sigma} \card{T_\sigma \setminus \varphi_j^{-1}(V_j)} = \sum_{j \in \sigma}\card{A_j \setminus V_j} \leq \sum_{j \in \sigma} \frac{\rho}{s}\card{G} \leq \rho \card{G}.
\end{align*}
Finally, if $i, j \in \sigma$, we have $(\varphi_i \compo \varphi_j^{-1})(Y_j) = (\varphi_i \compo \varphi_j^{-1})(A_j \setminus W_j) = \varphi_i(T_\sigma \setminus X_\sigma) = A_i \setminus W_i = Y_i$.
\end{proof}


\section{Sum-free self-inverse sets}
\label{sec:sum-free-self-inverse}

As a first application of \Cref{thm:polywrapper}, we obtain the following container theorem for self-inverse sum-free subsets of $\FF_q$.

\begin{theorem}
\label{thm:sum-free-container}
Let $\rL(q) \defined (\loglog q)^{1/5} / (\log q)^{1/10}$ and $0 < \eps < 1/4$.
There is a family $\cC$ of subsets of $\FF_q$ that satisfy:
\begin{enumerate}[label=(\roman*)]
    \item \label{item:container}
    for every $A \in \cF_\Sigma^\times(\FF_q)$, there $X \in \cC$ with $A \subseteq X$;
    \item \label{item:container-sparse}
    for every $X \in \cC$ we have $X = X^{-1}$ and $\card{X}/q \leq d_\Sigma(\FF_q)^2 + O_\eps\p{\rL(q)}$;
    \item \label{item:container-few}
    and $\log_2 \card{\cC} \ll q^{\eps}$.
\end{enumerate}
As a consequence, we have the following extremal result
\begin{equation}
\label{eq:finite-field-extremal}
    d_\Sigma^\times(\FF_q) \leq d_\Sigma(\FF_q)^2 + O\p{\rL(q)},
\end{equation}
and the following enumeration result
\begin{equation}
\label{eq:finite-field-counting}
    \frac{1}{q} \log_2 \card{\cF_\Sigma^\times(\FF_q)} \leq \frac{d_\Sigma(\FF_q)^2}{2} + O\p{\rL(q)}.
\end{equation}
We also have a supersaturation result: if $A \subseteq \FF_q$ satisfies $m_\Sigma(A) \leq \delta$ and $A = A^{-1}$, then
\begin{equation}
\label{eq:finite-field-supersaturation}
    \card{A}/q \leq d_\Sigma^\times(\FF_q) + O\p{ \min\set{\delta^{1/6}, \rL(q)} }.
\end{equation}
\end{theorem}
\begin{proof}
Let $\rho, \tau > 0$ be parameters satisfying $\tau \asymp \rL(q)^3$ and $\rho \asymp \rL(q)$.
We show that every set $A \in \cF_\Sigma^\times(\FF_q)$ can be approximated by an appropriate wrapper $W \subseteq \FF_q$.
Indeed, let $A \subseteq \FF_q^\times$ with $A = A^{-1}$ and $(\ind_A \conv \ind_A \conv \ind_{-A})(0) = 0$.
Apply \Cref{thm:polywrapper} with $A_1 = A_2 = A$, $A_3 = -A$, $\varphi_1(x) = \varphi_2(x) = x$, $\varphi_3(x) = -x$, partition $\Sigma = \set{\set{1,2,3}}$, and kernel $\kappa = q \ind_{\set{0}}$.
We obtain wrappers $W_1, W_2, W_3 \subseteq G$ and exceptional sets $Y_1, Y_2, Y_3 \subseteq G$ with $A \setminus Y_1 \subseteq W_1$, $A \setminus Y_2 \subseteq W_2$, and $-A \setminus Y_3 \subseteq W_3$.
Moreover, $W_2 = (\varphi_2 \compo \varphi_1^{-1})(W_1) = W_1$ and $W_3 = (\varphi_3 \compo \varphi_1^{-1})(W_1) = -W_1$, and hence if we set $W \defined W_1$ then we have $W = W_1 = W_2 = -W_3$.
Similarly, if $Y \defined Y_1$ then $Y = Y_1 = Y_2 = -Y_3$.

Thus, for every $A \in \cF_\Sigma^\times(\FF_q)$, we have shown that $A \setminus Y \subseteq W$, where $Y$ is a set with $\card{Y} \leq \rho q$, and $W$ is a $(\Lambda, K)$-wrapper with $ (\ind_W \conv \ind_W \conv \ind_{-W})(0) \leq \tau$ with 
\begin{equation*}
    K \ll \tau^{-1}\rho^{-2} \ll \rL(q)^{-5} = \frac{(\log q)^{1/2}}{\loglog q},
\end{equation*}
and the bound $\card{\Lambda} \leq D$, where we can take
\begin{equation*}
    D \ll \tau^{-2}\rho^{-4} \log(1/\tau \rho) \ll \rL(q)^{-10} \loglog q = \frac{\log q}{\loglog q}.
\end{equation*}
Since $A = A^{-1}$, we have that $A \setminus (Y \cup Y^{-1}) \subseteq W \cap W^{-1}$.
In other words, if we define the container $X \defined (W \cap W^{-1}) \cup Y \cup Y^{-1}$, we have $A \subseteq X$.
The family of containers we choose is the following:
\begin{equation*}
    \cC \defined \set[\big]{ (W \cap W^{-1}) \cup Y \cup Y^{-1} \st \text{$W$ is a $(\Lambda,K)$-wrapper with $\card{\Lambda} \leq D$ and $\card{Y} \leq \rho q$}}.
\end{equation*}

By definition, we have that $X = X^{-1}$ for $X \in \cC$ and we have already shown that for every $A \in \cF_\Sigma^\times(\FF_q)$, there is $X \in \cC$ with $A \subseteq X$, so \ref{item:container} is established.
Thus, we have
\begin{equation*}
    \card{A} \leq \card{X} \leq \card{W \cap W^{-1}} + \card{Y} + \card{Y^{-1}}.
\end{equation*}
To determine the cardinality of $W \cap W^{-1}$, write $w \defined \card{W}/q$ for the density of $W$ and note first that \Cref{prop:supersaturation-sumfree} gives $w \leq d_\Sigma(\FF_q) + O(\tau^{1/3})$, while \Cref{prop:kloosterman-mixing} gives
\begin{equation*}
    \abs[\big]{\card{W \cap W^{-1}} - w^2 q} \leq 2q^{1/2} \norm{\ind_W}_\wiener^2.
\end{equation*}
We will show later that $\norm{\ind_W}_\wiener \leq q^{\eps}$.
Given this, we have
\begin{align*}
    \card{X}/q &\leq w^2 + 2 \rho + 2 q^{-1/2} \norm{\ind_W}_\wiener^2 \\
    &\leq d_\Sigma(\FF_q)^2 + O(\tau^{1/3}) + 2 \rho + 2 q^{-1/4}
    \leq d_\Sigma(\FF_q)^2 + O_\eps(\rL(q)),
\end{align*}
as claimed in \ref{item:container-sparse}, which immediately gives \eqref{eq:finite-field-extremal} by fixing $\eps = 1/8$.

We now show that indeed $\norm{\ind_W}_\wiener \leq q^\eps$.
By \Cref{prop:wiener-wrapper}, the Wiener norm of $W$ satisfy $\norm{\ind_W}_\wiener \leq \p{K \cB(\FF_q)}^{\card{\Lambda}}$, while \Cref{cor:phi-G} gives $\cB(\FF_q) \leq 2 \log p$.
Therefore
\begin{align}
\label{eq:finite-field-wiener}
    \log \norm{\ind_W}_\wiener &\leq \card{\Lambda} \p[\big]{\log K + \log \cB(\FF_q)}
    \ll \tau^{-2} \rho^{-4} \log(1/\tau\rho) \p[\big]{\log(1/\tau \rho) + \loglog p} \ll \log q.
\end{align}
Therefore, $\tau$ and $\rho$ can be chosen such that $\log \norm{\ind_W}_\wiener \leq \eps \log q$, and thus $\norm{\ind_W}_\wiener \leq q^{\eps}$.

By \Cref{lemma:wrapper-count}, the number of $(\Lambda,K)$-wrappers in $G$ with $\card{\Lambda} \leq D$ is at most $2^{K^D + O(D \log \card{G})}$.
The number of subsets $Y \subseteq \FF_q$ with $\card{Y} \leq \rho q$ is $\sum_{i \leq \rho q}\binom{q}{i} \leq 2^{h(\rho)q}$ where $h$ is the binary entropy function $h(x) \defined - x \log_2(x) - (1-x)\log_2(1-x)$.
Together, this shows that
\begin{align*}
    \card{\cC} \leq 2^{K^D + O(D \log q)} \cdot 2^{h(\rho)}.
\end{align*}
Using the inequality $h(x) \leq x \log_2(4/x)$, we have $h(\rho) \ll \rho \log_2(4/\rho) \ll \rL(q) \loglog q$.
Furthermore, $D \log q \ll (\log q)^2$ and $K^D \leq q^{\eps}$ from \eqref{eq:finite-field-wiener}.
This gives
\begin{align*}
    \log_2 \card{\cC} \leq K^D + O(d \log q) + h(\rho) \ll q^\eps,
\end{align*}
which established \ref{item:container-few}.

For the enumeration result, note that every $X \in \cC$ has roughly $2^{\card{X}/2}$ subsets in $\cF_\Sigma^\times(\FF_q)$.
Indeed, fix a partition of $\FF_q^\times = I \sqcup I^{-1} \sqcup J$, where $J = \set{ x \in \FF_q^\times \st x^{-1} = x}$ and some $I$ with $\card{I} = (q-1-\card{J})/2$.
Note that $\card{J} = 1$ if $q$ is even and $\card{J} = 2$ if $q$ is odd.
Since $X = X^{-1}$, we have $\card{X \cap I} = \card{X}/2$ and every subset $A \subseteq X$ with $A = A^{-1}$ satisfy $A = (A \cap I) \sqcup (A \cap I)^{-1} \sqcup (A \cap J)$.
In other words, $A$ is determined solely by $A \cap (X \cap (I \cup J))$ and there are at most $2^{\card{X}/2 + 2}$ self-inverse subsets of $X$.
We can estimate then
\begin{equation*}
    \card{\cF_\Sigma^\times(\FF_q)} \leq \sum_{X \in \cC} 2^{\card{X}/2 + 2} \leq 4 \card{\cC} \; 2^{\max\set{\card{X} \st X \in \cC}/2}.
\end{equation*}
Using \ref{item:container-sparse}, we then have
\begin{equation*}
\label{}
    \frac{1}{q} \log_2 \card{\cF_\Sigma^\times(\FF_q)} \leq q^{-1+\eps} + \frac{d_\Sigma(\FF_q)^2}{2} + O_\eps\p{\rL(q)},
\end{equation*}
which gives \eqref{eq:finite-field-counting} by fixing $\eps = 1/8$.

Finally, to obtain the supersaturation result, note that we may assume that $\delta \gg \rL(q)^6$.
We apply the \Cref{thm:polywrapper} with different parameters.
Choose $\tau = \delta^{1/2}$ and $\rho = 16\, \delta^{1/6}$ so \eqref{eq:wrapper-parameters} is satisfied.
We get then that if $A \subseteq \FF_q$ is such that $A = A^{-1}$ and $(\ind_A \conv \ind_A \conv \ind_{-A})(0) \leq \delta$, then there is a $(\Lambda, K)$-wrapper $W$ with $(\ind_W \conv \ind_W \conv \ind_{-W})(0) \leq \delta^{1/2}$ and a exceptional set $Y$ with $\card{Y} \leq 16\, \delta^{1/6}$ with $A \subseteq W \cup Y$.
Thus we have $A \subseteq (W \cap W^{-1}) \cup Y \cup Y^{-1}$.
Using the supersaturation \Cref{prop:supersaturation-sumfree}, if $w \defined \card{W}/q$, we have $w \leq d_\Sigma(\FF_q) + O(\delta^{1/6})$ and thus $w^2 \leq d_\Sigma(\FF_q)^2 + O(\delta^{1/6})$.
Note also that $K \ll \delta^{-5/6}$ and $\card{\Lambda} \ll \delta^{-5/3}\log(1/\delta)$, so
\begin{align*}
    \log \norm{\ind_W}_\wiener \leq \card{\Lambda}\p[\big]{\log K + \loglog p)}
    \ll \delta^{-5/3} \log(1/\delta) \p[\big]{\log(1/\delta) + \loglog p} \ll \log q,
\end{align*}
and thus we may take $\norm{\ind_W}_\wiener \ll q^{1/8}$, which then implies $q^{1/2} \norm{\ind_W}_\wiener^2 \ll q^{3/4} \ll \delta^{1/6} q$.
Apply \Cref{prop:kloosterman-mixing} to obtain
\begin{equation}
\label{eq:supersatuartion-bound}
    \card{A} \leq \card{W \cap W^{-1}} + 2 \card{Y} \leq w^2 q + 2 q^{1/2} \norm{\ind_W}_\wiener^2 + 2\card{Y} \leq d_\Sigma(\FF_q)^2 q + O(\delta^{1/6}) q.
\end{equation}

While \eqref{eq:finite-field-extremal} only provides an upper bound, a corresponding on lower bound of
\begin{equation*}
    d_\Sigma^\times(\FF_q) \geq d_\Sigma(\FF_q)^2 - O(q^{-1/2}(\log p)^2)
\end{equation*}
can be obtained via the Kloosterman mixing~\Cref{prop:kloosterman-mixing} together with the construction from \Cref{lemma:sumfree-construction}.
Thus we can establish \eqref{eq:finite-field-supersaturation}, since\eqref{eq:supersatuartion-bound} turns to
\begin{equation*}
     \card{A} \leq d_\Sigma(\FF_q)^2 q + O(\delta^{1/6}) q \leq d_\Sigma^\times(\FF_q) q + O(\delta^{1/6}) q.
\end{equation*}
This also give a corresponding lower bound of $d_\Sigma(\FF_q)^2/2 - O(q^{-1/2}(\log p)^2)$ for \eqref{eq:finite-field-counting}.
\end{proof}


\section{Asymptotic independence}
\label{sec:indep}

In this section, we provide the upper bound four our general independence result \Cref{thm:general-independence-fq}.

\begin{theorem}
\label{thm:general-container}
Let $L_1, \dotsc, L_t \in \Lin_{\geq 3}$ and let $r_1, \dotsc, r_t$ be distinct non-zero integers.
Then for every $\eps > 0$, there is $q_0 = q_0(\eps) > 0$ such for all $q \geq q_0$ for which $L_1, \dotsc, L_t$ are non-degenerate in $\FF_q$ and $r_1, \dotsc, r_t$ are non-zero in $\FF_q$, then there is a family $\cC$ of subsets of $\FF_q$ that satisfy:
\begin{enumerate}[label=(\roman*)]
    \item \label{item:gen-container}
    for every $X \in \cF_{L_1, \dotsc, L_t}^{(r_1, \dotsc, r_t)}(\FF_q)$, there $C \in \cC$ with $X \subseteq C$;
    \item \label{item:gen-container-sparse}
    every $C \in \cC$ satisfy $\card{C}/q \leq d_{L_1}(\FF_q) \dotsb d_{L_s}(\FF_q) + o(1)$;
    \item \label{item:gen-container-few}
    and $\log_2 \card{\cC} \leq q^{o(1)}$.
\end{enumerate}
As a consequence, we have
\begin{equation}
\tag{\ref{eq:general-independence-fq}}
     d_{L_1, \dotsc, L_t}^{(r_1, \dotsc, r_t)}(\FF_q) = d_{L_1}(\FF_q) \dotsb d_{L_t}(\FF_q) + o(1) = q^{-1} \log_2 \card{\cF_{L_1, \dotsc, L_t}^{(r_1, \dotsc, r_t)}(\FF_q)}.
\end{equation}
\end{theorem}
\begin{proof}
Note that the lower bound in \eqref{eq:general-independence-fq} is established by \Cref{prop:general-construction} in \Cref{sec:constructions}.

We choose $\rho = \rho(q) > 0$ and $\tau = \tau(q) > 0$ satisfying $\rho(q) = o(1)$, $\tau(q) = o(1)$ and further constraints.
Let $K \defined c_s \tau^{-1} \rho^{1-s}$ and $D \defined c_s \tau^{-2} \rho^{2 - 2s} \log(1/\rho \tau)$ where $c_s$ and $c'_s$ are leading constants implied in \eqref{eq:wrapper-bound-k} and \eqref{eq:wrapper-bound-k} of \Cref{thm:polywrapper}.
Choose $\tau$ and $\rho$ such that $(K \cB(\FF_q))^D \leq q^{o(1)}$. 

Write $\cF \defined \cF_{L_1, \dotsc, L_t}^{(r_1, \dotsc, r_t)}(\FF_q)$ and let $X \in \cF$.
For each $1 \leq i \leq t$, we will build wrappers $W_i$ for the set $X^{r_i}$ as follows.
Since $X^{r_i}$ is $L_i$-free, we have $m_L(X^{r_1}) = 0$.
Write the linear form explicitly as $L_i(x_1, \dotsc, x_s) = a_1 x_1 + \dotsb + a_s x_s$ and let $A_1 \defined  a_1 \dilation X^{r_i}, \dotsc, A_s \defined a_s \dilation X^{r_i}$.
Note that $(\ind_{A_1} \conv \dotsb \conv \ind_{A_s})(0) = 0$.
Apply \Cref{thm:polywrapper} with partition $\Sigma = \set{\set{1,\dotsc,s}}$ and $\varphi_i(x) = a_i x$, which are automorphisms since $L_i$ is non-degenerate in $\FF_q$.
There is then a $(\Lambda_i, K)$-wrapper $W_i \subseteq \FF_q$ and an exceptional set $Y_i \subseteq \FF_q$ such that $X^{r_i} \setminus Y_i \subseteq W_i$, $\card{Y} \leq \rho q$, and
$m_L(W_i) = (\ind_{a_1 \dilation W_i} \conv \dotsb \conv \ind_{a_s \dilation W_i})(0) \leq \tau = o(1)$.
By \Cref{cor:green-removal}, we have $\card{W_i}/q \leq d_{L_i}(\FF_q) + o(1)$.
We also have that $\card{\Lambda_i} \leq D$, which by \Cref{prop:wiener-wrapper} gives
\begin{equation*}
    \norm{\ind_{W_i}}_\wiener \leq (K \cB(\FF_q))^D \leq q^{o(1)}.
\end{equation*}

Thus, for each $1 \leq i \leq t$, we have build wrappers $W_i$ of density $w_i \defined \card{W_i}/q$ and exceptional sets $Y_i$ of density $\card{Y_i}/q = o(1)$ satisfying $X^{r_i} \setminus Y_i \subseteq W_i$ and $\norm{\ind_{W_i}}_\wiener \leq q^{o(1)}$.
Consider then the set
\begin{equation*}
    Z \defined \set[\big]{x \in \FF_q^\times \st x^{r_1} \in W_1, \dotsc, x^{r_t} \in W_s },
\end{equation*}
and note that \Cref{prop:general-mixing} gives that 
\begin{equation*}
    \card{Z}/q \leq w_1 \dotsb w_t + c q^{-1/2} \norm{\ind_{W_1}}_\wiener \dotsb \norm{\ind_{W_t}}_\wiener \leq w_1 \dotsb w_t + o(1).
\end{equation*}
However, note that $X$ satisfies
\begin{align*}
    \card{X} &\leq \card[\big]{\set{x \in \FF_q^\times \st x^{r_1} \in X^{r_1}, \dotsc, x^{r_t} \in X^{r_t} }} 
    \leq \card{Z} + \sum_{1 \leq i \leq t} \card[\big]{\set{x \in \FF_q^\times \st x^{r_i} \in Y_i}} \\
    &= \card{Z} + \sum_{1 \leq i \leq t} \gcd(\abs{r_i},q-1) \card{Y_i}
    = \p[\big]{w_1 \dotsb w_t + o(1) }q,
\end{align*}
where we have used that the homomorphism $x \mapsto x^r$ in $\FF_q^\times$ has kernel of size $g = \gcd(r,q-1)$, so $\card{\set{x \in \FF_q^\times \st x^r \in Y}} \leq g \card{Y}$, attained when $Y$ has no more than one element in each fibre.
In any case, this shows that for every $X \in \cF$, we have
\begin{equation*}
    \card{X}/q \leq w_1 \dotsb w_t + o(1) \leq d_{L_1}(\FF_q) \dotsb d_{L_t}(\FF_q) + o(1).
\end{equation*}
Together with the lower bound in \Cref{prop:general-construction}, this establishes the density part of \eqref{eq:general-independence-fq}.

Thus, if we define the following family of containers
\begin{equation*}
    \cC \defined \set[\bigg]{ \set[\big]{x \in \FF_q^\times \st x^{r_1} \in W_1, \dotsc, x^{r_t} \in W_t} \cup Y  \st \begin{array}{c}
    \text{$W_1, \dotsc, W_t$ are $(\Lambda,K)$-wrappers} \\
    \text{with $\card{\Lambda} \leq D$ and $\card{Y} \leq \rho' q$}
  \end{array} },
\end{equation*}
where $\rho' = \rho(\abs{r_1} + \dotsb + \abs{r_t})$, then we have shown that for every $X \in \cF$, there is $C \in \cC$ with $X \subseteq C$, so \ref{item:gen-container} is established.
Moreover, we indeed have also showed that $\card{C}/q \leq d_{L_1}(\FF_q) \dotsb d_{L_t}(\FF_q) + o(1)$, which gives \ref{item:gen-container-sparse}.
We now enumerate the containers.
By \Cref{lemma:wrapper-count}, the number of $(\Lambda,K)$-wrappers in $\FF_q$ with $\card{\Lambda} \leq D$ is at most $2^{K^D + O(D \log q)}$ and he number of subsets $Y \subseteq \FF_q$ with $\card{Y} \leq \rho' q$ is $\sum_{i \leq \rho q}\binom{q}{i} \leq 2^{h(\rho')q}$.
Together, this shows \ref{item:gen-container-few} since
\begin{align*}
    \log_2 \card{\cC} \leq s K^D + O(s D \log q) + h(\rho) \leq q^{o(1)}.
\end{align*}

Finally, we enumerate $\cF$ and obtain the upper bound in \eqref{eq:general-independence-fq}.
Note that $\card{\cF} \leq \card{\cC} 2^{M}$ where $M = \max\set{\card{C} \st C \in C}$.
Therefore 
\begin{equation*}
    \log_2 \card{\cF} \leq \log_2 \card{C} + M \leq d_{L_1}(\FF_q) \dotsb d_{L_t}(\FF_q) + o(1). \qedhere
\end{equation*}
\end{proof}

\Cref{thm:general-independence-fq} is finally proven.

In the next proposition, we record for completeness an enumeration of the $L$-free subsets of $\FF_q$ whenever $L$ has only two variables.

\begin{proposition}
\label{prop:two-variables}
Let $L \in \Lin_2(\FF_q)$ and write $L(x,y) = ax + by$ for non-zero $a,b \in \FF_q$.
Let $\gamma \defined - a/b$ and let $k$ be the order of $\gamma$ in $\FF_q^\times$.
Then we have $d_L(\FF_q) = (1-1/q)\floor{k/2}/k$ and $\card{\cF_L(\FF_q)} = L_k^{(q-1)/k}$, where $L_k$ is the $k$-th Lucas number.
\end{proposition}
\begin{proof}
Note that an $L$-free subset $A \subseteq \FF_q$ must not have $0$.
Moreover, for $x \in \FF_q^\times$, if $x \in A$ then $\gamma x \notin A$ and $\gamma^{-1} x \notin A$.
Let $G$ be the Cayley graph with vertex set $\FF_q^\times$ and multiplicative generators $\set{\gamma, \gamma^{-1}}$.
An $L$-free set corresponds precisely to an independent subset in $G$.
Note that $G$ is a collection of $q/k$ disjoint cycles of length $k$, each corresponding to a different orbits of the action $x \mapsto \gamma x$.
A cycle $C_k$ of length $k$ has $L_k$ independent subsets, so $G$ has $L_k^{(q-1)/k}$ independent subsets.
The largest independent set in $C_k$ has size $\floor{k/2}$, so $d_L(\FF_q) = (1-q)\floor{k/2}/k$.
\end{proof}

For a fixed linear form $L(x,y) = ax + by$, \Cref{prop:two-variables} implies that $d_L(\FF_p)$ and $(1/p)\log_2\card{\cF_L(\FF_p)}$ do not converge as $p \to \infty$ unless $\gamma \in \set{-1,1}$.
These limits exists, for instance, along a subsequence in which $\gamma$ is a primitive root in $\FF_q^\times$, and we would have $d_L(\FF_q) = 1/2 + o(1)$ and $(1/q)\log_2\card{\cF_L(\FF_q)} = \log_2((1 + \sqrt{5})/2) + o(1)$.
In any case, we cannot have $(1/q)\log_2\card{\cF_L(\FF_q)}  = d_L(\FF_q) + o(1)$ (that is, \eqref{eq:card-family-L-Fq}) unless $a = -b$ and $d_L(\FF_q) = 0$.

If $L(x,y) = ax + by$ with $a = b = 0$, then in fact $L \in \Lin_0$ and $\cF_L(\FF_q) = \emptyset$ so $d_L(\FF_q) = \card{\cF_L(\FF_q)} = 0$.
If only one of $a$ or $b$ is non-zero, then $\cF_L(\FF_q) = \set{A \subseteq \FF \st 0 \notin A}$ so $d_L(\FF_q) = 1 - 1/q$ and $\card{\cF_L(\FF_q)} = 2^{d_L(\FF_q) q} = 2^{q-1}$.


\section{Product-free sets}
\label{sec:product-free}

In this section, we prove \Cref{thm:general-independence-fq-product}, and as a corollary, we obtain \Cref{thm:sum-product}.
The proof shares a lot with the proof of \Cref{thm:general-independence-fq}.

\begin{theorem}
\label{thm:general-container-product}
Let $L_1, \dotsc, L_t \in \Lin_{\geq 3}$ and let $r_1, \dotsc, r_t$ be distinct non-zero integers.
Then for every $\eps > 0$, there is $q_0 = q_0(\eps) > 0$ such for all $q \geq q_0$ for which $L_1, \dotsc, L_t$ are non-degenerate in $\FF_q$ and $r_1, \dotsc, r_t$ are non-zero in $\FF_q$, then there is a family $\cC$ of subsets of $\FF_q$ that satisfy:
\begin{enumerate}[label=(\roman*)]
    \item \label{item:gen-prod-container}
    for every $X \in \cF_{L_1, \dotsc, L_t, \Pi}^{(r_1, \dotsc, r_t)}(\FF_q)$, there $C \in \cC$ with $X \subseteq C$;
    \item \label{item:gen-prod-container-sparse}
    every $C \in \cC$ satisfy $\card{C}/q \leq d_{L_1}(\FF_q) \dotsb d_{L_s}(\FF_q) d_\Pi(\FF_q^\times) + o(1)$;
    \item \label{item:gen-prod-container-few}
    and $\log_2 \card{\cC} \leq q^{o(1)}$.
\end{enumerate}
As a consequence, we have
\begin{equation}
\tag{\ref{eq:general-independence-fq-product}}
     d_{L_1, \dotsc, L_t, \Pi}^{(r_1, \dotsc, r_t)}(\FF_q) = d_{L_1}(\FF_q) \dotsb d_{L_t}(\FF_q)  d_\Pi(\FF_q^\times) + o(1) = q^{-1} \log_2 \card{\cF_{L_1, \dotsc, L_t, \Pi}^{(r_1, \dotsc, r_t)}(\FF_q)}.
\end{equation}
\end{theorem}
\begin{proof}
The lower bound in \eqref{eq:general-independence-fq-product} is established by \Cref{prop:general-construction-product}.
We now prove the upper bound by following the proof of \Cref{thm:general-independence-fq} directly.
Indeed, let $\rho = \rho(q) > 0$, $\tau = \tau(q) > 0$,  $K \defined c_s \tau^{-1} \rho^{1-s}$, and $D \defined c_s \tau^{-2} \rho^{2 - 2s} \log(1/\rho \tau)$ satisfy $\rho(q) = o(1)$, $\tau(q) = o(1)$, $(K \cB(\FF_q))^D \leq q^{o(1)}$ so we can apply \Cref{thm:polywrapper}.

Write $\cF \defined \cF_{L_1, \dotsc, L_t, \Pi}^{(r_1, \dotsc, r_t)}(\FF_q)$ and let $X \in \cF$.
Following \Cref{thm:general-independence-fq}, for each $1 \leq i \leq t$, we have build wrappers $W_i$ of density $w_i \defined \card{W_i}/q \leq d_{L_i}(\FF_q) + o(1)$ and exceptional sets $Y_i$ of density $\card{Y_i}/q = o(1)$ satisfying $X^{r_i} \setminus Y_i \subseteq W_i$ and $\norm{\ind_{W_i}}_{\wiener(\FF_q)} \leq q^{o(1)}$.

We now note that $X \subseteq \FF_q^\times$ and that $X$ is product-free.
Therefore, in the commutative group $\FF_q^\times$, we have $(\ind_X \conv \ind_X \conv \ind_{-X})(0) = 0$.
We apply \Cref{thm:polywrapper} with $G = \FF_q^\times$, $\Sigma = \set{\set{1,2,3}}$, $\varphi_1(x) = \varphi_2(x) = x$, $\varphi_3(x) = - x$, and we obtain a wrapper $W \subseteq \FF_q^\times$ and an exceptional set $Y$ with $X \setminus Y \subseteq W$, with $\norm{\ind_W}_\wiener \leq q^{(1)}$ and $m_\Sigma(W) = o(1)$, which by \Cref{prop:supersaturation-sumfree}, implies that $w \defined \card{W}/q \leq d_\Sigma(\FF_q^\times) + o(1)$.

Now that we have all the wrappers at our disposal, consider the sets
\begin{align*}
    Z &\defined \set[\big]{x \in \FF_q^\times \st x^{r_1} \in W_1, \dotsc, x^{r_t} \in W_t, x \in W }, \\
    J &\defined Y \cup \set[\big]{x \in \FF_q^\times \st \text{ $x^{r_i} \in Y_i$ for some $1 \leq i \leq t$}  },
\end{align*}
and note that $\card{J} \leq \rho' q$ where $\rho' = \rho(1 + \abs{r_1} + \dotsb + \abs{r_t})$.
Crucially, we have also that $X \in Z \cup J$.
In other words, every $X \in \cF$ is contained in one set of the following family of containers
\begin{equation*}
    \cC \defined \set[\bigg]{ \set[\big]{x \in \FF_q^\times \st x^{r_1} \in W_1, \dotsc, x^{r_t} \in W_t, x \in W} \cup J  \st \begin{array}{c}
    \text{$W_1, \dotsc, W_t$, $W$ are $(\Lambda,K)$-wrappers} \\
    \text{with $\card{\Lambda} \leq D$ and $\card{J} \leq \rho' q$}
  \end{array} },
\end{equation*}
and \ref{item:gen-prod-container} is established.
Note that \Cref{prop:general-mixing} gives that for every $C \in \cC$, we have \ref{item:gen-prod-container-sparse} since
\begin{align*}
    \card{C}/q &\leq w_1 \dotsb w_t w  + c q^{-1/2} \norm{\ind_{W_1}}_{\wiener(\FF_q)} \dotsb \norm{\ind_{W_t}}_{\wiener(\FF_q)} \norm{\FF_q^\times}_{\wiener(\FF_q^\times)} \\
    &\leq w_1 \dotsb w_t w + o(1) \leq d_{L_1}(\FF_q) \dotsb d_{L_s}(\FF_q) d_\Pi(\FF_q^\times) + o(1).
\end{align*}
As in \Cref{thm:general-independence-fq}, we have $\log_2 \card{\cC} \leq q^{o(1)}$, and therefore \ref{item:gen-prod-container-few} is established.
Finally, we enumerate $\cF$ and obtain the upper bound in \eqref{eq:general-independence-fq-product}.
Write $M = \max\set{\card{C} \st C \in C}$ and again
\begin{equation*}
    \log_2 \card{\cF} \leq \log_2 \card{C} + M \leq d_{L_1}(\FF_q) \dotsb d_{L_s}(\FF_q)d_\Pi(\FF_q^\times) + o(1).
\end{equation*}
This finished the proof of \eqref{eq:general-independence-fq-product}.
\end{proof}

For completeness, we finish by determining that the set of limiting densities in $d_\Pi(\FF_{2^d})$.

\begin{proposition}
\label{prop:produc-liminfsup}
The largest density of a product-free subset of $\FF_{2^d}$ are such that
\begin{equation}
\tag{\ref{eq:liminfsup}}
    \frac{1}{3} = \liminf_{d \to \infty} d_\Pi(\FF_{2^d}) < \limsup_{d \to \infty} d_\Pi(\FF_{2^d}) = \frac{2}{5}.
\end{equation}
\end{proposition}
\begin{proof}
Let $N = 2^d -1$ and recall $d_\Pi(\FF_{2^d}) = (1 - 1/2^d) d_\Sigma(\ZZ/N\ZZ)$.
Since $N$ is odd, every prime divisor of $N$ that is congruent to $2 \mod{3}$ is at least $5$, the formula for $d_{\Sigma}(\ZZ /N \ZZ)$ gives
\begin{equation*}
    \frac{1}{3} - \frac{1}{3N} \leq d_{\Sigma}(\ZZ / N\ZZ) \leq \frac{2}{5}.
\end{equation*}
We can find subsequences for the values of $d$ that give the desired limits.
First note that if $d = 4k$, then $5 \divides N$ and hence we attain the limit of $2/5$.
For the inferior limit, choose $d$ to be prime.
Indeed, if $p$ is a prime that divisor of $2^d -1$, then Fermat's little theorem gives that $d \divides p - 1$, and thus $p > d$. 
Then, the formula for $d_{\Sigma}(\ZZ / N \ZZ)$ gives $d_\Sigma(\ZZ/N\ZZ) \leq 1/3 + 1/3d$.
Therefore, along primes the limit tends to $1/3$ as desired. 
\end{proof}


\printbibliography

\end{document}